\documentclass[10pt,a4paper]{article}

\usepackage[utf8]{inputenc}
\usepackage{amsmath,amsfonts,amssymb,amsthm,hyperref}
\usepackage{a4wide}
\usepackage{nomencl}
\makenomenclature
\usepackage{nicefrac}
\usepackage{graphicx}
\usepackage{caption}
\usepackage{subcaption}
\graphicspath{{figures/}}
\DeclareGraphicsExtensions{.pdf}

\newcommand{\dt}{\Delta t}
\newcommand{\dx}{\Delta X}
\newcommand{\dy}{\Delta Y}
\newcommand{\hf}{\frac{1}{2}}
\newcommand{\lx}{x_{j(i)}-x_i}
\newcommand{\ly}{y_{j(i)}-y_i}
\newcommand{\vx}{u_{j(i)}-u_i}
\newcommand{\vy}{v_{j(i)}-v_i}
\newcommand{\Lag}{\mathsf{L}}
\newcommand{\uref}{U_{\text{ref}}}
\newcommand{\vref}{V_{\text{ref}}}
\newcommand{\norm}[1]{\left\lVert#1\right\rVert}
\newcommand\abs[1]{\left|#1\right|}

\DeclareMathOperator*{\argmin}{arg\,min}

\newtheorem{remark}{Remark}

\begin{document}

\title{Macroscopic modeling of multi-lane motorways using a two-dimensional second-order model of traffic flow}
\date{\today}
\author{Michael Herty \\
		{\small\it Institut f\"{u}r Geometrie und Praktische Mathematik (IGPM)} \\
		{\small\it RWTH Aachen University} \\
		{\small\it Templergraben 55, 52062 Aachen, Germany\vspace{5mm}} \\ 
	Salissou Moutari \\
		{\small\it School of Mathematics and Physics} \\
		{\small\it Mathematical Sciences Research Centre} \\
        {\small\it Queens' University, Belfast} \\      
		{\small\it BT7 1NN, Northern Ireland, United Kingdom \vspace{5mm}} \\ 
	Giuseppe Visconti \\
		{\small\it Institut f\"{u}r Geometrie und Praktische Mathematik (IGPM)} \\
		{\small\it RWTH Aachen University} \\
		{\small\it Templergraben 55, 52062 Aachen, Germany}
	   }

\maketitle

\begin{abstract}
Lane changing is one of the most common maneuvers on motorways. Although, macroscopic traffic models are well known for their suitability to describe fast moving crowded traffic, most of these models are generally developed in a one dimensional framework, henceforth lane changing behavior is either not explicitly modeled or explicitly forbidden. In this paper, we propose a macroscopic model, which accounts for lane-changing behavior on motorway, based on a two-dimensional extension of the Aw and Rascle~\cite{aw2000SIAP} and Zhang~\cite{ZhangMacro} macroscopic model for traffic flow.  Under conditions, when lane changing maneuvers are no longer possible, the model ``relaxes'' to the one-dimensional Aw-Rascle-Zhang model. Following the same approach as in \cite{aw2002SIAP}, we derive the two-dimensional macroscopic model through  scaling of  time discretization of a microscopic follow-the-leader model with driving direction. We provide a detailed analysis of the space-time discretization of the proposed macroscopic model as well as an approximation of the solution to the associated Riemann problem. Furthermore, we illustrate some promising features of the proposed model through some numerical experiments.
\end{abstract}



\section{Introduction} \label{sec:Introduction}

Current macroscopic models for multi-lane traffic on motorways generally couple a multi-lane or multi-class one-dimensional first-order model of traffic flow (LWR model)~\cite{lighthill1955PRSL,richards1956OR} with some lane-changing rules in order to capture traffic dynamics stemming from the lane-changing maneuvers~\cite{Daganzo2002a,Greenbergetal2003,HollandandWoods1997,LavalandDaganzo2006,Michalopoulosetal1984}. In our opinion, this approach suffers from the following aspects:
\begin{enumerate}
\item A major assumption of the LWR model is that for a given traffic density all the drivers adopt the same velocity. However, this assumption is not always valid in practice. Macroscopic second-order models, e.g.~\cite{aw2000SIAP,ZhangMacro}, attempted to address this limitation by introducing an additional variable, which in~\cite{FanHertySeibold} was interpreted as the ``relative'' velocity of a specific (class) of drivers. The advection of the ``relative'' velocity with the actual velocity enable  to describe the reaction of drivers to traffic conditions ahead. Therefore, second-order models can be viewed as an extension of the LWR such that, for a given traffic density, drivers can react differently, by adopting a wide range of speed~\cite{FanHertySeibold,Lebacque03,LebacqueGSOM}. In~\cite{aw2002SIAP}, the derivation of the second-order macroscopic models, presented in ~\cite{aw2000SIAP,ZhangMacro}, was proposed. More recently, Di Francesco et al \cite{DiFrancescoFagioliRosini} proposed a ``follow-the-leader'' approximation of the Aw-Rascle-Zhang model in a multi-population framework together with a detailed discussion on the analytical properties of the model.
\item Multi-lane models, based on the coupling of one-dimensional equations for each lane, describe the lane-changing behavior through some interaction terms (namely sources on the right-hand sides of the equations) derived from empirical interaction rates, see e.g.~\cite{klar1999SIAP-1, klar1999SIAP-2, KlarWegener2000ab, HoogendoornBovy, Helbing, Zhuetal}. The interaction is typically assumed to be proportional to local density on both the current and the target lane. A fluid dynamics model describing the cumulative density on all lanes is proposed in~\cite{ChetverushkinChurbanovaFurmanovTrapeznikova2010, ChetverushkinChurbanovaSukhinovaTrapeznikova2008, SukhinovaTrapeznikovaChetverushkinChurbanova2009}, where a two-dimensional system of balance laws is obtained by analogy with the quasi-gas-dynamics (QGD) theory. In the aforementioned studies, the two-dimensional modeling of traffic dynamics assumed that vehicles move to lanes with faster speed or lower density, and the evolution for the lateral velocity is described as proportional to the local density and the mean speed along the road. A major shortcoming of such approach is the estimation of the interaction rate from data as well as the corresponding high number of parameters. In order to address these limitations,in this study we consider a different approach, namely we treat lanes as continuum avoiding the need to prescribe heuristically the dynamics of the flow of vehicles across lanes.

\end{enumerate}

In a recent work, Herty and Visconti~\cite{HmVg2017} used microscopic data, namely vehicle trajectories collected on a German highway, to derive a two-dimensional first-order macroscopic model. The study presented in this paper used the same experimental microscopic data to derive a macroscopic model for multi-lane traffic on highways based on a two-dimensional extension of the Aw-Rascle-Zhang model~\cite{aw2000SIAP,ZhangMacro}. By revisiting the analysis proposed in~\cite{aw2002SIAP}, we show that the semi-discretization of the two-dimensional Aw-Rascle-Zhang model can be viewed as the limit of a multi-lane ``follow-the-leader'' model. These results enable the two-dimensional macroscopic model to capture traffic dynamics caused by the lane-changing behavior in spite a coarse scale traffic flow description. Through some numerical simulations, we highlight the relationship between the both models as well as their ability to reproduce classical traffic situations. The numerical results presented in this study are obtained using a classical first order finite volume scheme since our main objective was to show the relationship between the microscopic and the macroscopic description of traffic flow. However, suitable numerical strategies should be used to approximate the macroscopic equations e.g. the scheme proposed in~\cite{GoatinChalons} which enables to reduce oscillations nearby contact discontinuities.

\medskip

Note that in~\cite{AlNasurThesis, AlNasurKachroo}, a derivation of a two-dimensional macroscopic crowd model from a microscopic model was proposed. In contrast to these studies, here we provide a detailed analysis of the approximation of the solution to the associated Riemann problem. The approach adopted in the current paper was to approximate the solution to the two-dimensional Riemann problem through a coupling of two-half Riemann problems in one dimension and then use the transition from one lane to another, subject to some admissibility conditions, in order to capture the two-dimensional aspect of the problem. We provide a detailed discussion, of the main scenarios under consideration in this study, as well as the suitable Riemann solvers at the interfaces of lanes, which describe the lateral dynamics caused by lane-changing maneuvers.  A distinctive feature of the proposed model is its double-sided behavior. Under conditions, when lane-changing maneuvers are no longer possible, e.g. the traffic is congested on adjacent lanes, the model ``relaxes'' to the one-dimensional Aw-Rascle-Zhang model~\cite{aw2000SIAP,ZhangMacro}, whereas in free flow conditions on adjacent lanes, the model captures dynamics caused by car movement between adjacent lanes.

%Besides a preprint \cite{HmVg2017} on 
%a first-order macroscopic model for two-dimensional traffic flow we are not aware of further publications on truely two-dimensional macroscopic models. The derivation of the macroscopic model from precise microscopic dynamics has been conducted for spatially one-dimensional models already in the paper \cite{aw2002SIAP}. To some extend we follow this derivation, however, this has to be modified to treat the lateral dynamics. 

\medskip

The remaining part of this paper is organized as follows. 
Section~\ref{sec:preliminary} presents a brief overview of the one-dimensional Aw-Rascle-Zhang~\cite{aw2000SIAP,ZhangMacro} model, reviewing some of its mathematical properties relevant to the current study. In particular, we revisit the derivation of the Aw-Rascle-Zhang model from the microscopic ``follow-the-leader'' model using a different approach. Section~\ref{sec:2DARmodel} and Section~\ref{sec:LagrToEul} introduce the proposed 2D macroscopic model, outline and discuss sufficient conditions for the derivation of the model through some scaling of the time discretization of a microscopic ``follow-the-leader'' model with driving direction. A detailed analysis of the two-dimensional macroscopic model as well as an approximation of the solution to the associated Riemann problem are provided in Section~\ref{sec:2DMacroModel.properties}. In Section~\ref{sec:numericalSimulations}, we discuss the space-time discretization of the proposed macroscopic and we present some numerical simulations highlighting some distinctive features of the proposed model. Finally, Section~\ref{sec:conclusion} concludes the paper and outlines some directions for further research.

\section{Preliminary discussions} \label{sec:preliminary}

In this section, we will start by reviewing the study presented in~\cite{aw2002SIAP}, where the authors show a connection between the classical microscopic ``follow-the-leader'' (FTL) model and the Aw and Rascle~\cite{aw2000SIAP} and Zhang~\cite{ZhangMacro} (ARZ) macroscopic model. More precisely, they prove that the FTL model can be viewed as semi-discretization of the ARZ model in Lagrangian coordinates. However, here our aim is to prove the connection between the two models using a different approach, namely in the case of two-dimensional models.

In contrast with the analysis \cite{aw2002SIAP}, here we consider a two-dimensional FTL model first and derive the corresponding macroscopic limiting equations using suitable coordinate transformations. In order to illustrate the idea we briefly recall the conservative form of the second-order traffic flow model \cite{aw2000SIAP, ZhangMacro} without relaxation term:
\begin{equation} \label{eq:1D:ARZ}
	\begin{cases}
		\partial_t \rho(t,x) + \partial_x (\rho u)(t,x) =0,\\
		\partial_t (\rho w)(t,x) + \partial_x (\rho u w)(t,x)=0,
	\end{cases}
\end{equation}
where $w(t,x)=u(t,x)+P$ and $P=P(\rho)$ is the so called ``traffic pressure''.

The corresponding one--dimensional microscopic model is based on the following arguments. The movement of particles is described according to $\dot{x}_i =u_i$, where $x_i=x_i(t)$ and $u_i=u_i(t)$ are the time-dependent position and speed, respectively, of vehicle $i$. If not specified, the quantities are assumed at time $t$. In Lagrangian coordinate the particles
move along with the trajectory. The local density at time $t$ is defined as 
\begin{equation} \label{eq:1Ddensity}
	\rho_i = \frac{\dx}{\lx},
\end{equation}
where $\dx$ is the typical length of a vehicle and $\tau_i=1/\rho_i$ describes the dynamics of the ``relative distance''. The index $j(i)$ identifies the interacting vehicle for a car $i$. Let's assume, without loss in generality, that cars are ordered such that $x_i < x_{i+1}$, $\forall\,i$. Then, under this hypothesis,  we  have $j(i)=i+1$. Therefore, the distance between the $i^\text{th}$ and the $(i+1)^\text{th}$ car evolves as follows:
\begin{equation} \label{EvolutionMicroTau}
	\tau_i(t+\dt)-\tau_i(t) = \frac{x_{i+1}(t+\dt) - x_{i+1}(t) - ( x_i(t+\dt) - x_i(t))}{\dx}.
\end{equation}
Dividing both sides of (\ref{EvolutionMicroTau}) by $\dt$ and taking the limit yields the following equation: 
\begin{equation} \label{eq:1D:EvolutionMicroTau}
	\dot{\tau}_i = \frac{u_{i+1}-u_i}{\dx}.
\end{equation}
We define a macroscopic velocity by 
$
	u_{i}(t) =: u^{\Lag}(t,X_i),
$
where $X_i=X(t,x_i)$ is the cumulative car mass up to the car labeled $i$ at time $t$. By definition 
\begin{equation} \label{eq:1Dcarmass}
	X(t,x) = \int^x \rho(t,\xi) d\xi.
\end{equation}
Notice that we are assuming that time is not influenced by the transformation from the discrete dynamics. In a discrete model with $N$ cars, we define the density around car $i$ as follows:
$$
	\rho(t,x)= \frac{1}{\dx} \sum^N_{k=1} \delta(x-x_k(t)).
$$
The total mass is then $N / \dx$. Therefore, we obtain:
$$
	X(t,x_i) = \frac{1}{\dx} \left(\sum^{i-1}_{l=1} 1 + \hf\right).
$$
Since the cars are ordered such that $x_i<x_{i+1}$,then for each fixed time $t$, the map $s \mapsto X(t,s)$ is a monotone function for a given density. Therefore, there is a one-to-one map from $i$ to $X(t, x_i).$ Hence, we adopt the notation 
$X_i = X(t,x_i).$ Afterwards, we can extend the values $u_i$ to a function $u^{\Lag}$ such that, at each time $t$, $u_i(t) = u^{\Lag}(t,X_i).$ We proceed similarly for $\tau_i.$ 
\noindent
Finally, from equation \eqref{eq:1D:EvolutionMicroTau} we obtain the following:
\begin{equation*}
	\partial_t \tau^{\Lag}(t,X_i) = \frac{u^{\Lag}(t,X_{i+1})-u^{\Lag}(t,X_{i})}{\dx},
\end{equation*}
which is a finite volume semi-discretization of the PDE
\begin{equation} \label{eq:1DLagrConsLaw}
	\partial_t \tau^{\Lag} - \partial_X u^{\Lag} = 0,
\end{equation}
under the following assumptions:
\begin{description}
	\item[Ansatz 1:] we have one car per cell and the distance between the centers of two cells is precisely $\dx$, i.e. $\dx$ is the grid space.
	\item[Ansatz 2:] $\tau^{\Lag}$ and $u^{\Lag}$ are distributed as piecewise constant in each cell $\left[X_{i-\frac12},X_{i+\frac12}\right]$.
\end{description}
Note that the first assumption means that we are considering the macroscopic limit where the number of vehicles goes to infinity while the length of cars shrinks to zero. The second assumption implies that we are considering a first-order finite--volume scheme and without serious restrictions.

Now, let's to consider the equation for the acceleration in the ``follow-the-leader'' model without relaxation, given by:
$$
	\dot{u}_i = U_\text{ref} \dx^\gamma \frac{v_{i+1}-v_i}{(x_{i+1}-xi)^{\gamma+1}}.
$$
Let $w_i = u_i + P(\tau_i)$ where $P(\tau_i)$ is a function defined as follows:
$$
	P(\tau_i) = \begin{cases}
		\displaystyle{\frac{U_\text{ref}}{\gamma\tau_i^\gamma}}, & \gamma>0\\
		\displaystyle{-U_\text{ref} \ln(\tau_i}), & \gamma=0.
	\end{cases}
$$
Then, straightforward computations show that $\dot{w}_i = \dot{u}_i + P^\prime(\tau_i)\dot{\tau}_i = 0$, which can represent the second equation of the particle model. Again, we can identify $w_i(t) = w^\Lag(t,X_i)$ and thus we get
\begin{equation} \label{eq:1DLagrSpeed}
	\partial_t w^\Lag = 0.
\end{equation}

As proved in~\cite{aw2002SIAP}, equation~\eqref{eq:1DLagrConsLaw} and equation~\eqref{eq:1DLagrSpeed} give the Lagrangian version of the ARZ model~\eqref{eq:1D:ARZ}. Thus, the microscopic FTL model can be derived through a semi-discretization of the macroscopic ARZ model in Lagrangian coordinates.

\begin{figure}[t!]
	\centering
	\includegraphics[width=\textwidth]{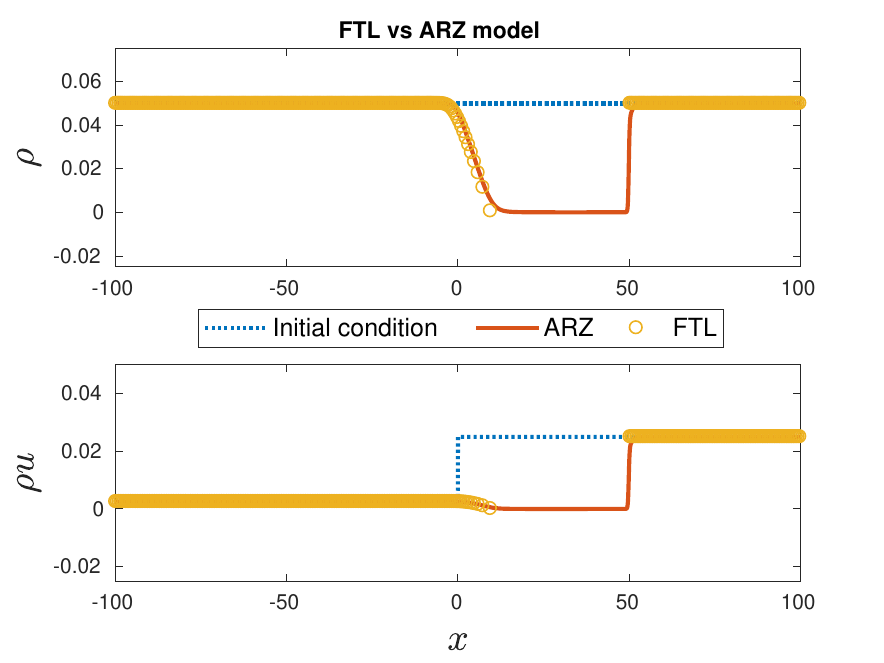}
	\caption{Simulation of the 1D ARZ model (red solid line) and of the 1D FTL model (orange circles). The top panel shows the density, while the bottom panel shows the flux at initial time (dotted line) and at final time.\label{fig:1D}}
\end{figure}

In Figure~\ref{fig:1D}, we replicate the simulation results presented in~\cite{aw2002SIAP} in order to give a numerical evidence of the relationship between the two models. The blue dotted line is the initial condition for the density $\rho$ and the flux $\rho u$, while the red solid line and the orange circles give the solution of the macroscopic and microscopic model, respectively. We refer the reader to~\cite{aw2002SIAP} for further details on the simulation parameters.

\section{Derivation of a two-dimension spatial extension of the ARZ model} \label{sec:2DARmodel}

In this section, we will derive the extension of the ARZ model to the case of two space dimensions. To achieve this, we will first introduce a generalization of the FTL model with dynamics including lane changing and then we will use the similar arguments to Section~\ref{sec:preliminary} to derive the macroscopic model.

A two-dimensional microscopic model requires the evolution, in time, of the positions along and also across the road section. Let $x_i=x_i(t)$ and $y_i=y_i(t)$ be the time-dependent positions. Then the evolution of the positions is given by:
$$
	\dot{x}_i=u_i, \quad \dot{y}_i=v_i,
$$
for any vehicle $i=1,\dots,N$. The speed $u_i$ is supposed to be non-negative (travel in $x$-direction), while the speed $v_i$ can be either positive or negative (travel in $y$-direction).

\begin{figure}[t!]
	\centering
	\includegraphics[width=\textwidth]{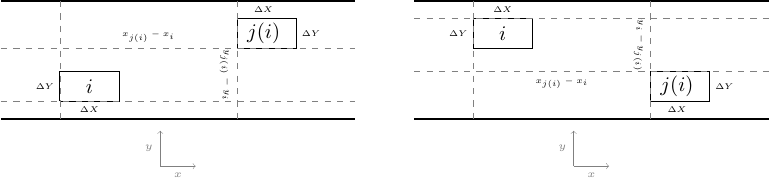}
	\caption{Different position between the test vehicle $i$ and the interacting vehicle~$j(i)$. Left: the $y$-distance is $y_{j(i)}-y_i>0$. Right: the $y$-distance is $y_i-y_{j(i)}>0$.\label{fig:2D:Density}}
\end{figure}

In contrast with the one-dimensional case, here we assume that there is no particular ordering imposed right now, and we simply label the cars. The space occupied by each car is given by $\dx \dy$, where $\dx$ and $\dy$ are the typical length and the typical width of a vehicle, respectively. Generalizing the definition of the one-dimensional case, the density around vehicle $i$ becomes
\begin{equation} \label{eq:2Ddensity}
	\rho_i = \frac{\dx\dy}{\left(\lx\right)\abs{\ly}},
\end{equation}
where the vehicle, interacting with vehicle $i$, is labeled $j(i)$. Since here we do not assume any order on the labeling, in general $j(i) \neq i+1$. The absolute value in \eqref{eq:2Ddensity} is considered in order to take into account both  situations depicted in Figure~\ref{fig:2D:Density}, depending on the relative position of vehicle $j(i)$ with respect to vehicle $i$, and thus to guarantee the positivity of the density. Using the coordinate system in Figure~\ref{fig:2D:Density}, we assume that a vehicle $i$ moves towards the left side of the road if $v_i>0$ and towards the right side if $v_i<0$.

\begin{remark}
	Let $\rho_i^{\text{1D}}$ and $\rho_i^{\text{2D}}$ be the local density in the one-dimensional (see \eqref{eq:1Ddensity}) and in the two-dimensional case (see \eqref{eq:2Ddensity}), respectively. Observe that
	$$
		\rho_i^{\text{2D}} \to \rho_i^{\text{1D}}, \quad \text{as} \quad \dy, \abs{\ly} \to 0^+
	$$
	only if we assume that
	\begin{equation} \label{eq:LimitHyp}
	\lim_{\substack{\dy \to 0^+ \\ (\ly) \to 0^+ }} \frac{\dy}{\abs{\ly}} = 1.
	\end{equation}
\end{remark}

\begin{figure}[t!]
	\centering
	\includegraphics[width=\textwidth]{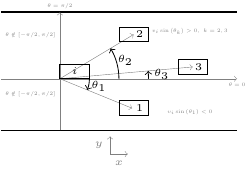}
	\caption{Choice of the interacting car in the case $v_i>0$. The interacting vehicle will be car $2$, namely the nearest vehicle in the driving direction of vehicle $i$.\label{fig:2D:ChoiceOfj}}
\end{figure}

Let us define $Q_i :=(x_i,y_i)$, $i=1,\dots,N$, the vector of the microscopic positions of vehicle $i$. Then, the interacting vehicle $j(i)$ is determined by the following map:
\begin{equation} \label{eq:FieldCar}
	i \mapsto j(i) = \argmin_{\substack{h=1,\dots,N \\ v_i \sin\theta_h > 0 \\ \theta_h\in\left[-\frac{\pi}{2},\frac{\pi}{2}\right]}} \norm{Q_h-Q_i}_2.
\end{equation}
The motivation behind this choice is as follows (see also Figure~\ref{fig:2D:ChoiceOfj}). Assume that each test vehicle $i$ defines a coordinate system in which the origin is its right rear corner if $v_i\geq 0$ and its left rear corner if $v_i<0$. We are indeed dividing the road in four areas. Let $\theta_h$ be the angle between the $x$-axis (in the car coordinate system) and the position vector $Q_h$ of vehicle $h$. Then, the request $\theta_h\in\left[-\frac{\pi}{2},\frac{\pi}{2}\right]$ enables to consider only cars being in front of vehicle $i$. On the other hand, the request $v_i\sin\theta_h>0$ allows to consider only cars in the driving direction of vehicle $i$. Among all these vehicles we choose the nearest one. Therefore, \eqref{eq:FieldCar} can be rewritten as follows:
$$
i \mapsto j(i) = \argmin_{\substack{h=1,\dots,N \\ v_i (y_h - y_i)>0 \\ x_h > x_i}} \norm{Q_h-Q_i}_2.
$$

\noindent
The specific volume for car $i$ is defined as $\tau_i : = 1/\rho_i$ and thus we have 
$$
	\tau_i = \frac{ \left(\lx\right)\abs{\ly}}{\dx \dy}.
$$
From elementary geometry, we obtain, as in the one-dimensional case, the following relationship:
$$
	\tau_i(t+\dt) - \tau_i(t) = \frac{ \left(\lx\right)\abs{\ly} (t+\dt) }{\dx \dy} - \frac{ \left(\lx\right)\abs{\ly} (t) }{\dx \dy}.
$$
Now, we add and subtract $\frac{ \left(\lx\right)(t) \abs{\ly} (t+\dt) }{\dx \dy}$, and tend $\dt$ to zero to obtain the final dynamics: 
% todo{Notice that the absolute value allows to take into account the ``direction of the wind''.}
\begin{equation} \label{eq:2D:EvolutionMicroTau}
	\dot{\tau}_i = \frac{ \left(\vx\right)\abs{\ly} }{ \dx \dy} + \frac{ \left(\vy\right) \abs{\ly} \left(\lx\right) }{ \dy\dx \left(\ly\right) }.
\end{equation}

As in the one-dimensional case, we aim to identify $\tau_i$ with a function out of the discrete dynamics in order to find the corresponding conservation law. Then, we introduce 
$$
X = \int^x \rho(t,\xi,\bar{\eta}) d\xi, \quad Y = \int^y \rho(t,\bar{\xi},\eta) d\eta
$$
for some values $\bar{\eta}$, $\bar{\xi}$. $X$ and $Y$ are the cumulative car mass up to the car labeled $i$ in $x$-direction (and projected on $y=\bar{\eta}$) and in $y$-direction (and projected on $x=\bar{\xi}$), respectively.
 
Again, we would like to relate the label of the car to a corresponding mass. 
The discrete measure corresponding to $N$ cars is now given by:
$$
	\rho(t,x,y) = \frac{1}{\dx\dy} \sum^N_{i=1} \delta(x-x_i(t)) \delta(y-y_i(t)),
$$
with total mass $N / \dx / \dy.$ 

\medskip

In order to obtain a one-to-one relation between the label $i$ and the pair of indices $(k,\ell)$ corresponding to the cumulative masses $X_k$ and $Y_\ell$, we proceed as follows. We view the cars as points in the 2D domain. Then, we fix an arbitrary value $\bar{\eta}$ in $y$-direction and project all cars towards this line, i.e., from $(x_i,y_i)$ to $(x_i,\bar{\eta})$. 
Then, the projected density is given by: 
$$
	\rho(t,x,\bar{\eta}) = \frac{1}{\dx\dy} \sum_{i=1}^N \delta(x-x_i(t)).
$$
Afterwards, computing $X$ on the projected density gives
$$
	X(t,x_i,\bar{\eta}) = \frac{1}{\dx\dy} \left( \sum_{h \mid x_h\in(-\infty,x_i]} 1 + \hf \right) = \frac{1}{\dx\dy} \left( k+\hf \right),
$$
for some $k$.

\medskip

Similarly, we fix $\bar{\xi}$ in $x$ direction and project all cars towards this line, i.e., from $(x_i,y_i)$ to $(\bar{\xi},y_i)$. We obtain 
$$
	Y(t,\bar{\xi},y_i) = \frac{1}{\dx\dy} \left( \sum_{h \mid x_h\in(-\infty,y_i]} 1 + \hf \right) = \frac{1}{\dx\dy} \left( \ell+\hf \right)
$$
for some $\ell$.

\medskip

Note that $k$ and $\ell$ may well be different according to the position of the vehicles. Moreover, for each fixed time $t$, the definition does not depend on the choice of $\bar{\eta}$ and $\bar{\xi}$.

Again, in the projected densities, the maps $s\mapsto X(t,s,y)$ and $s\mapsto Y(t,x,s)$ are monotone functions. Hence, for those
quantities, at least for short time, there is a one-to-one correspondence between car $i$ and the pair $(k,\ell)$, respectively. Therefore, we can identify 
$$
	\tau_i =\tau^{\Lag}(t,X_k,Y_\ell),
$$
and similarly for the speeds
$$
	u_{i}=u^{\Lag}(t,X_{k},Y_{\ell}), \quad v_{i}=v^{\Lag}(t,X_{k},Y_{\ell}).
$$

Since we  are interested in the limit for many cars, it is natural to assume that there are enough cars in the neighboring cells. We therefore assume that:
\begin{description}
	\item[Ansatz 3:] $j(i) \mapsto (k+1,\ell+1)$ if $v_i\geq 0$ \ and \ $j(i) \mapsto (k+1,\ell-1)$ if $v_i< 0$.
\end{description}

Now, observe that, since $X$ and $Y$ are monotone, we can write 
$$
	x_i = X^{-1}\left( \frac{1}{\dx\dy} \left(k+\hf\right) \right), \quad y_i = Y^{-1}\left( \frac{1}{\dx\dy} \left(\ell+\hf\right) \right).
$$
Assume that $v_i$ is positive, then by straightforward computation using the fact that the map $Y$ is locally invertible, we obtain that: 
$$
	y_{j(i)} - y_i  = Y^{-1}\left( \frac{1}{\dx\dy} \left(\ell+1+\hf\right) \right) - Y^{-1}\left( \frac{1}{\dx\dy} \left(\ell+\hf\right) \right) = \dy.
$$

%In fact, we compute
%\begin{align*}
%	Y^{-1}\left( \frac{1}{\dx\dy} \left(\ell+1+\hf\right) \right) - Y^{-1}\left( \frac{1}{\dx\dy} \left(\ell+\hf\right) %\right) = \dy \\
%	\Leftrightarrow \frac{1}{\dx\dy} \left(\ell+1+\hf\right) = Y \left( \dy + Y^{-1}\left( \frac{1}{\dx\dy} \left(\ell+%\hf\right) \right) \right) \\
%	\Leftrightarrow \frac{1}{\dx\dy} \left(\ell+1+\hf\right) = Y \left( \dy + y_i \right) \\
%	\Leftrightarrow \frac{1}{\dx\dy} \left(\ell+1+\hf\right) = \left( \int^{y_i} + \int_{y_i}^{y_i+\dy} \right) \rho(t,%\bar{\xi},\eta) d\eta  \\
%	\Leftrightarrow \frac{1}{\dx\dy} \left(\ell+1+\hf\right) = \frac{1}{\dx\dy} \left(\ell+\hf\right) + \int_{y_i}^{y_i+\dy} \rho(t,\bar{\xi},\eta) d\eta \\
%	\Leftrightarrow \frac{1}{\dx\dy} \left(\ell+1+\hf\right) = \frac{1}{\dx\dy} \left(\ell+\hf\right) + \frac{1}{2\dx\dy} + \frac{1}{2\dx\dy},
%\end{align*} 
%where the last is true provided that we have $y_{j(i)} = y_i+\dy$ (which is reasonable in the macroscopic limit) and that %there are no other $\delta$'s between $y_i$ and $y_i+\dy$. 

Similar computations hold for the case $v_i<0$ and in the $x$-direction. Therefore, equation~\eqref{eq:2D:EvolutionMicroTau} rewrites as follows:
\begin{equation} \label{eq:2D:FinalEvolutionMicroTau}
	\dot{\tau}_i = \frac{ \left(\vx\right) }{\dx} + \frac{ \left(\vy\right) \abs{\ly} }{ \dy \left(\ly\right) }.
\end{equation}
Now, we show that~\eqref{eq:2D:FinalEvolutionMicroTau} can be seen as a suitable semi-discretization in space of~\eqref{eq:2DLagrConsLaw}.
\begin{equation} \label{eq:2DLagrConsLaw}
	\partial_t \tau^{\Lag} - \partial_X u^{\Lag} - \partial_Y v^{\Lag} = 0.
\end{equation}
%which is the straightforward generalization of~\eqref{eq:1DLagrConsLaw}. In other words, we prove %that~\eqref{eq:2D:FinalEvolutionMicroTau} can be seen as a suitable semi-discretization in space %of~\eqref{eq:2DLagrConsLaw}.

\noindent
We consider a uniform Cartesian grid
$$
	\Omega_{k,\ell} = \left[X_{k-\hf},X_{k+\hf}\right] \times \left[Y_{\ell-\hf},Y_{\ell+\hf}\right]	
$$
in which $X_{k+\hf}=X_k+\dx$ and $Y_{\ell+\hf}=Y_\ell+\dx$ and we define the volume average as follows:
$$
	\overline{\tau^\Lag}_{k,\ell} (t) = \frac{1}{\dx\dy} \iint_{\Omega_{k,\ell}} \tau^\Lag(t,\xi,\eta) \mathrm{d}\xi\mathrm{d}\eta.
$$
Integrating the conservation law over each control volume and dividing by the volume of $\Omega_{k,\ell}$, we obtain the finite volume formulation of ~\eqref{eq:2DLagrConsLaw}:
\begin{equation} \label{eq:FV2DLagrConsLaw}
\begin{aligned}
	\frac{\mathrm{d}}{\mathrm{d}t} \overline{\tau^\Lag}_{k,\ell} = & \frac{1}{\dx\dy} \int_{Y_{\ell-\hf}}^{Y_{\ell+\hf}} \left( u^\Lag(t,X_{k+\hf},\eta)-u^\Lag(t,X_{k-\hf},\eta) \right) d\eta \\ & + \frac{1}{\dx\dy} \int_{X_{k-\hf}}^{X_{k+\hf}} \left( v^\Lag(t,\xi,Y_{\ell+\hf})-v^\Lag(t,\xi,Y_{\ell-\hf}) \right) d\xi.
\end{aligned}
\end{equation}
We consider a first-order finite--volume scheme. This implies that $\tau^\Lag$, $u^\Lag$ and $v^\Lag$ are given as piecewise constant on each of the above patches. In this interpretation, the mass of the car is uniformly distributed on the path $\dx \times \dy$ and thus, for the specific volume we have
$$
	\overline{\tau^\Lag}_{k,\ell}(t) = \tau^\Lag(t,X_k,Y_\ell).
$$

\noindent
Since the velocity along the road is non-negative, then for the first term in the right-hand side of~\eqref{eq:FV2DLagrConsLaw}, an Upwind flux would be appropriate and in this case $u^\Lag(t,X_{k+\hf},\eta) \approx u^\Lag(t,X_{k+1},\eta)$. Then, we write 
$$
	\frac{1}{\dx\dy}  \int_{Y_{\ell-\hf}}^{Y_{\ell+\hf}} \left( u^\Lag(t,X_{k+1},\eta)-u^\Lag(t,X_{k},\eta) \right) d\eta. 
$$
Since $u^\Lag$ is constant on the patch, we evaluate it at any point $\eta \in \left[Y_{\ell-\hf},Y_{\ell+\hf}\right]$. We choose
to evaluate it now at different points in $Y$ direction, so that for the first term in the right-hand side of~\eqref{eq:FV2DLagrConsLaw} we finally get:
$$
	\frac{1}{\dx}  \Big( u^\Lag(t,X_{k+1},Y_{\ell+1})-u^\Lag(t,X_{k},Y_{\ell})  \Big)
$$

\noindent
For the second term in the right-hand side of~\eqref{eq:FV2DLagrConsLaw}, we notice that the speed across the lanes can be either positive or negative. We again use an Upwind flux. However, if the velocity is positive, then the approximation $v^\Lag(t,\xi,Y_{\ell+\hf}) \approx v^\Lag(t,\eta,Y_{\ell+1})$ is still appropriate. On the other hand, if the velocity is negative, then we use the following approximation $v^\Lag(t,\xi,Y_{\ell+\hf}) \approx v^\Lag(t,\eta,Y_{\ell})$. Thus, we obtain:
$$
	\begin{cases}
		\displaystyle{\frac{1}{\dx\dy}  \int_{X_{k-\hf}}^{X_{k+\hf}} \left( v^\Lag(t,\xi,Y_{\ell+1})-v^\Lag(t,\xi,Y_{\ell}) \right) d\xi}, & \text{if \ $v^\Lag(t,X_k,Y_\ell)\geq 0$}, \\[0.5cm]
		\displaystyle{\frac{1}{\dx\dy}  \int_{X_{k-\hf}}^{X_{k+\hf}} \left( v^\Lag(t,\xi,Y_{\ell})-v^\Lag(t,\xi,Y_{\ell-1}) \right) d\xi}, & \text{if \ $v^\Lag(t,X_k,Y_\ell)< 0$}.
	\end{cases}
$$
Since $v^\Lag$ is also constant on the patch, we can evaluate it at any point $\xi \in \left[X_{k-\hf},X_{k+\hf}\right]$. We choose to evaluate it at different points in $X$ direction, so that we can finally obtain
$$
\begin{cases}
	\displaystyle{\frac{1}{\dy}  \Big( v^\Lag(t,X_{k+1},Y_{\ell+1})-v^\Lag(t,X_k,Y_{\ell}) \Big)}, & \text{if \ $v^\Lag(t,X_k,Y_\ell)\geq 0$} \\[0.5cm]
	\displaystyle{\frac{1}{\dy}  \Big( v^\Lag(t,X_k,Y_{\ell})-v^\Lag(t,X_{k-1},Y_{\ell-1}) \Big) }, & \text{if \ $v^\Lag(t,X_k,Y_\ell)< 0$}.
\end{cases}
$$

\noindent
Putting together all the terms, from~\eqref{eq:FV2DLagrConsLaw} we have
\begin{align*}
	\frac{\mathrm{d}}{\mathrm{d}t} \tau^\Lag_{k,\ell} = & \frac{1}{\dx}  \Big( u^\Lag_{k+1,\ell+1}(t)-u^\Lag_{k,\ell}(t)  \Big) + \frac{v^+}{\dy}  \Big( v^\Lag_{k+1,\ell+1}(t)-v^\Lag_{k,\ell}(t) \Big) + \frac{v^-}{\dy}  \Big( v^\Lag_{k,\ell}(t)-v^\Lag_{k-1,\ell-1}(t) \Big),
\end{align*}
where $v^+=\max\left(0,\frac{\abs{v^\Lag_{k,\ell}(t)}}{v^\Lag_{k,\ell}(t)}\right)$ and $v^-=\min\left(0,\frac{\abs{v^\Lag_{k,\ell}(t)}}{v^\Lag_{k,\ell}(t)}\right)$. Thus, we \eqref{eq:2D:FinalEvolutionMicroTau} is  a first-order finite volume semi-discretization of the conservation law~\eqref{eq:2DLagrConsLaw}.

Now, we study what happens for the acceleration. Actually, since a two-dimensional FTL model has not already been introduced in literature, we do not have the evolution equation for the acceleration in the two-dimensional case at hand. Thus, we first  need to derive the equations for $\dot{u}_i$ and $\dot{v}_i$.

We define two quantities
$$
	w_i = u_i + P_1(\tau_i), \quad \sigma_i = v_i + P_2(\tau_i),
$$
which can be seen as desired speeds of vehicle $i$ in the $x$- and in the $y$-direction, respectively. The quantities $P_1$ and $P_2$ play the role of the ``traffic pressure'' and they are functions of the local density. However, observe that $P_1$ and $P_2$ should be homogeneous to the velocity and for this reason we take $P_1\neq P_2$, assuming that there are two different reference velocities in the two different directions of the flow. We define
\begin{equation}\label{functionsP}
	P_1(\tau_i)=\begin{cases}
		\displaystyle{\frac{\uref}{\gamma_1 \tau_i^{\gamma_1}}}, & \text{if \ $\gamma_1 > 0$}\\[0.3cm]
		-\uref\ln(\tau_i), & \text{if \ $\gamma_1 = 0$}
	\end{cases}, \quad
	P_2(\tau_i)=\begin{cases}
	\displaystyle{\frac{\vref}{\gamma_2 \tau_i^{\gamma_2}}}, & \text{if \ $\gamma_2 > 0$}\\[0.3cm]
	-\vref\ln(\tau_i), & \text{if \ $\gamma_2 = 0$}
	\end{cases},
\end{equation}
and as in the one-dimensional model without relaxation, we require that the desired speeds are constant during the time evolution, that is
$$
	\dot{w}_i = \dot{\sigma}_i = 0.
$$
We therefore obtain
$$
	\dot{u}_i = -\dot{\tau}_i P_1^\prime(\tau_i), \quad \dot{v}_i = -\dot{\tau}_i P_2^\prime(\tau_i),
$$
and computing the derivatives of $P_1$ and $P_2$ with respect to $\tau_i$, we can finally obtain the following equations for the evolution of the microscopic accelerations:
%\todo[inline]{I do not like the following formulation because if we imagine that, e.g., $\abs{\ly}=0$ then the velocity of the vehicle $i$ tends to increase even if $\lx \gg 0$, i.e. the two vehicles are distant but in the same lane. This is strictly dependent on the choice of the local density $\rho_i$ which is, however, the proper way to generalize the one-dimensional density.}
\begin{align*}
	\dot{u}_i & = C_1 \left( \frac{\vx}{\left(\lx\right)\Delta A^{\gamma_1}}  + \frac{\vy}{\left(\ly\right)\Delta A^{\gamma_1}}\right),	\\
	\dot{v}_i & = C_2 \left( \frac{\vx}{\left(\lx\right)\Delta A^{\gamma_2}}  + \frac{\vy}{\left(\ly\right)\Delta A^{\gamma_2}}\right),
\end{align*}
where
$$
C_1=\uref\dx^{\gamma_1}\dy^{\gamma_1}, \quad C_2=\vref\dx^{\gamma_2}\dy^{\gamma_2}, \quad \Delta A = \left(\lx\right)\abs{\ly}.
$$

\begin{remark}
	Note that the above equations for the microscopic accelerations are consistent with
	\begin{align*}
		\dot{u}_i & = \uref\dx^{\gamma_1} \frac{\vx}{\left(\lx\right)^{\gamma_1+1}},	\\
		\dot{v}_i & = 0,
	\end{align*}
	which are the acceleration equations in the one-dimensional case in the limit $\vref \to 0$ and under hypothesis~\eqref{eq:LimitHyp}.
\end{remark}

\begin{remark}
	In the macroscopic limit, i.e. when the number of cars increases whereas $x_{j(i)}-x_i=\dx$ and $\abs{y_{j(i)}-y_i}=\dy$, we simply get
	\begin{align*}
		\dot{u}_i & = \uref \left( \frac{\vx}{\left(\lx\right)}  + \frac{\vy}{\left(\ly\right)}\right),	\\
		\dot{v}_i & = \vref \left( \frac{\vx}{\left(\lx\right)}  + \frac{\vy}{\left(\ly\right)}\right).
	\end{align*}
\end{remark}

Finally, the two-dimensional ``follow-the-leader'' microscopic model is given by the following equatios:
\begin{equation} \label{eq:2DFTL}
\begin{aligned}
	\dot{x}_i &= u_i,\\
	\dot{y}_i &= v_i,\\
	\dot{u}_i & = C_1 \left( \frac{\vx}{\left(\lx\right)\Delta A^{\gamma_1}}  + \frac{\vy}{\left(\ly\right)\Delta A^{\gamma_1}}\right),	\\
	\dot{v}_i & = C_2 \left( \frac{\vx}{\left(\lx\right)\Delta A^{\gamma_2}}  + \frac{\vy}{\left(\ly\right)\Delta A^{\gamma_2}}\right),
\end{aligned}
\end{equation}
which can be, altogether, rewritten in the form
\begin{equation} \label{eq:2DFTLLagr}
\begin{aligned}
	\dot{\tau}_i & = \frac{ \left(\vx\right) }{\dx} + \frac{ \left(\vy\right) \abs{\ly} }{ \dy \left(\ly\right) }, \\ 
	\dot{w}_i &= 0, \\
	\dot{\sigma}_i & = 0.
\end{aligned}
\end{equation}

We have already seen that the equation for the time evolution of the specific volume (i.e. the first equation in~\eqref{eq:2DFTLLagr}) is a first-order finite volume semi-discretization of the conservation law in Lagrangian coordinates~\eqref{eq:2DLagrConsLaw}. Since there exists a one-to-one map $i\mapsto(k,\ell)$, we may extend the values $w_i$ and $\sigma_i$ to two functions $w^\Lag(t,X,Y)$ and $\sigma^\Lag(t,X,Y)$ such that
$$
	w_i(t) = w^\Lag(t,X_k,Y_\ell), \quad \sigma_i(t) = \sigma^\Lag(t,X_k,Y_\ell),
$$
where $w^\Lag$ and $\sigma^\Lag$ are defined as follows:
$$
	w^\Lag(t,X,Y)=u^\Lag(t,X,Y)+P_1\left(\tau^\Lag(t,X,Y)\right), \quad \sigma^\Lag(t,X,Y)=v^\Lag(t,X,Y)+P_2\left(\tau^\Lag(t,X,Y)\right).
$$
Assuming a first-order scheme and thus assuming that $w^\Lag$ and $\sigma^\Lag$ are constant on the grid, $\dot{w}_i=0$ and $\dot{\sigma}_i=0$ can obviously be viewed as the results of a semi-discretization of
\begin{equation} \label{eq:2Dmac}
\begin{aligned}
	\partial_t w^\Lag = 0, \quad \partial_t \sigma^\Lag=0.
\end{aligned}
\end{equation}
Finally, the above equation~\eqref{eq:2Dmac} and~\eqref{eq:2DLagrConsLaw} define the following system of macroscopic equations in Lagrangian coordinates:
\begin{align}
	\partial_t \tau^\Lag = \partial_X u^\Lag + \partial_Y v^\Lag, \quad
	\partial_t w^\Lag  = 0, \label{eq:2D:LagrangianMacroModel}, \quad
	\partial_t \sigma^\Lag  = 0.
\end{align}

\begin{remark}[Validation of the 2D particle model]
	In~\cite{HmVg2017} we used some experimental data collected on a German highway in order to derive a two-dimensional first-order macroscopic model. The microscopic trajectories of vehicles were used to compute macroscopic quantities related to traffic flow and consequently to fundamental diagrams. The latter were used to define a closure for the macroscopic equation. Here, we use the same microscopic data introduced in~\cite{HmVg2017} to validate the two-dimensional microscopic model~\eqref{eq:2DFTL}. To achieve this, we fix an initial time and using experimental data we get the initial positions and speeds of vehicles on the road. We estimate the characteristics of the interacting vehicle for each car using~\eqref{eq:FieldCar} and then we evolve the trajectories using~\eqref{eq:2DFTL}. The right-most vehicle in $x$-direction is used as ghost car and its trajectory is updated at each time using the real trajectory. In the left panel of Figure~\ref{fig:ValidationMicro} we show the real trajectories (red line) provided by experimental data and the computed trajectories with the 2D microscopic model (black symbols), after $2$ seconds. In the right panel of Figure~\ref{fig:ValidationMicro} we show the error between the real and the computed trajectory for each car, using the 2-norm distance.
	\begin{figure}
		\centering
		\includegraphics[width=0.49\textwidth]{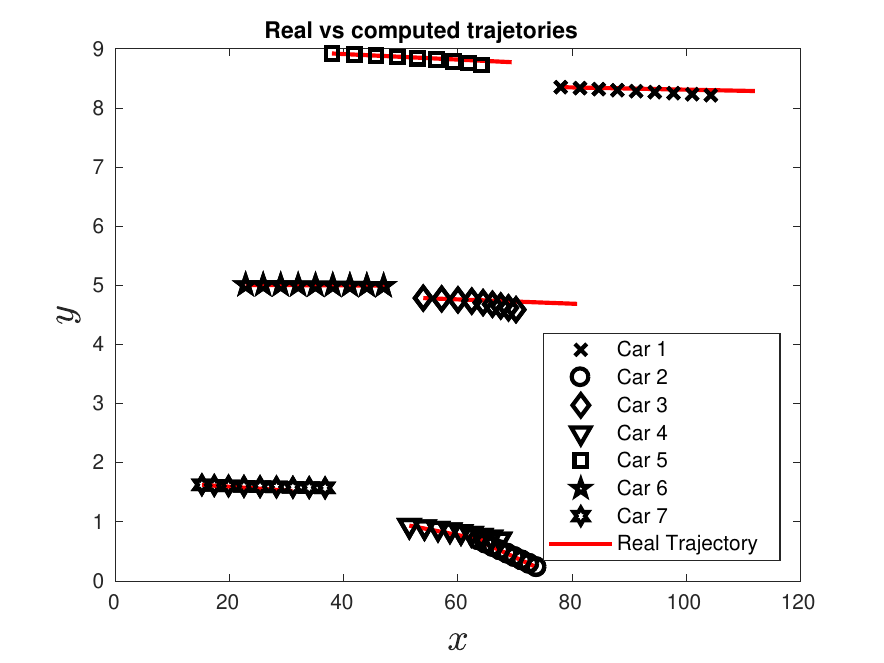}
		\includegraphics[width=0.49\textwidth]{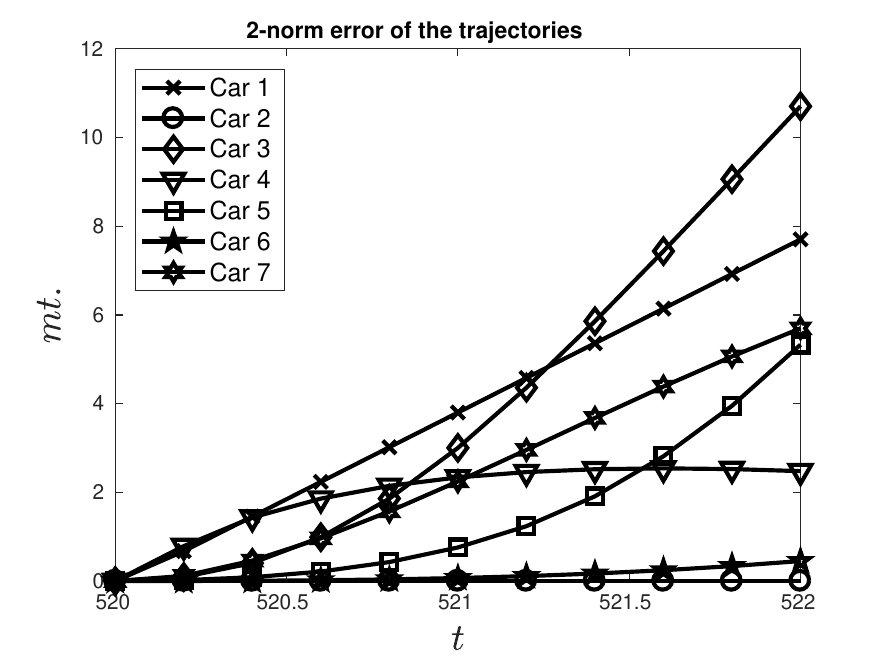}
		\caption{Left: comparison between the real trajectories provided by experimental data in~\cite{HmVg2017} and the computed trajectories using the two-dimensional microscopic model~\eqref{eq:2DFTL}. Right: evolution of the 2-norm error in time for each car.\label{fig:ValidationMicro}}
	\end{figure}
\end{remark}

\section{From Lagrangian to Eulerian coordinates} \label{sec:LagrToEul}

Macroscopic equations \eqref{eq:2D:LagrangianMacroModel} are written with respect to the ``mass'' coordinates (or Lagrangian) $X$ and $Y$, where the variables $X(t,x,y)$ and $Y(t,x,y)$ denote the total mass of vehicles up to point $x$ and up to point $y$, for $y=\overline{y}$ and $x=\overline{x}$ fixed, respectively.

We reformulate the model in Eulerian coordinates $x$ and $y$ by setting $t(T,X,Y)=T$ and 
$$
	x(t,X,\bar{Y}) = \int^X \tau^\Lag(t,\xi,\bar{Y}) \mathrm{d}\xi, \quad y(t,\bar{X},Y) = \int^Y \tau^\Lag(t,\bar{X},\eta) \mathrm{d}\eta.
$$
From this definition, and using the continuity equation in Lagrangian coordinates,  we obtain
\begin{gather*}
	\partial_X x = \partial_Y y = \tau^\Lag, \quad \partial_t x = u^\Lag, \quad \partial_t y = v^\Lag, \quad \partial_Y x=\partial_X y=0. 
\end{gather*}
The Jacobian of the coordinate transformation is given by 
\[
	J=\frac{\partial(T,X,Y)}{\partial(t,x,y)}=
%	\begin{bmatrix}
%	\partial_T t & \partial_X t & \partial_Y t \\[1ex]
%	\partial_T x & \partial_X x & \partial_Y x \\[1ex]
%	\partial_T y & \partial_X y & \partial_Y y
%	\end{bmatrix}
%	=
	\begin{bmatrix}
	1 & 0 & 0 \\[1ex]
	u^\Lag & \tau^\Lag & 0 \\[1ex]
	v^\Lag & 0 & \tau^\Lag
	\end{bmatrix},
\]
%Therefore the conversion is given by
%\[
%	\begin{pmatrix}
%	\displaystyle{\frac{\partial}{\partial T}}\\[1.5ex]
%	\displaystyle{\frac{\partial}{\partial X}}\\[1.5ex]
%	\displaystyle{\frac{\partial}{\partial Y}}
%	\end{pmatrix}
%	= J^\mathsf{T} \begin{pmatrix}
%	\displaystyle{\frac{\partial}{\partial t}}\\[1.5ex]
%	\displaystyle{\frac{\partial}{\partial x}}\\[1.5ex]
%	\displaystyle{\frac{\partial}{\partial y}}
%	\end{pmatrix}
%\]
and therefore 
%\[
%\begin{cases}
%	\displaystyle{\frac{\partial}{\partial T} = \frac{\partial}{\partial t} \frac{\partial t}{\partial T} + \frac{\partial}%{\partial x} \frac{\partial x}{\partial T} + \frac{\partial}{\partial y} \frac{\partial y}{\partial T}}\\[1.5ex]
%	\displaystyle{\frac{\partial}{\partial X} = \frac{\partial}{\partial t} \frac{\partial t}{\partial X} + \frac{\partial}%{\partial x} \frac{\partial x}{\partial X} + \frac{\partial}{\partial y} \frac{\partial y}{\partial X}}\\[1.5ex]
%	\displaystyle{\frac{\partial}{\partial Y} = \frac{\partial}{\partial t} \frac{\partial t}{\partial Y} + \frac{\partial}%{\partial x} \frac{\partial x}{\partial Y} + \frac{\partial}{\partial y} \frac{\partial y}{\partial Y}}
%\end{cases}
%	\Longrightarrow
%\begin{cases}
%	\displaystyle{\frac{\partial}{\partial T} = \frac{\partial}{\partial t} + u^\Lag \frac{\partial}{\partial x} + v^\Lag %\frac{\partial}{\partial y}}\\[1.5ex]
%	\displaystyle{\frac{\partial}{\partial X} = \tau^\Lag \frac{\partial}{\partial x}}\\[1.5ex]
%	\displaystyle{\frac{\partial}{\partial Y} = \tau^\Lag \frac{\partial}{\partial y}}
%\end{cases}.
%\]
%
%Let $\tau$, $u$ and $v$ be the Eulerian quantities related to the specific volume $\tau^\Lag$ and the speeds $u^\Lag$, %$v^\Lag$, that is
%\begin{align*}
%	\tau^\Lag(T,X,Y)&=\tau\big(t,x(t,X,Y),y(t,X,Y)\big),\\
%	u^\Lag(T,X,Y)&=u\big(t,x(t,X,Y),y(t,X,Y)\big),\\
%	v^\Lag(T,X,Y)&=v\big(t,x(t,X,Y),y(t,X,Y)\big).
%\end{align*}
%Using
the continuity equation in Lagrangian coordinates~\eqref{eq:2DLagrConsLaw} yields the corresponding equation 
in Eulerian variables
\[
	\frac{\partial \tau}{\partial t} + u \frac{\partial \tau}{\partial x} + v \frac{\partial \tau}{\partial y} = \tau \frac{\partial u}{\partial x} + \tau \frac{\partial v}{\partial y}.
\]
Recalling that $\tau=1/\rho$, we get the following equation:
%\[
%	-\frac{1}{\rho^2}\frac{\partial \rho}{\partial t} -\frac{1}{\rho^2}u \frac{\partial \rho}{\partial x} -\frac{1}{\rho^2} %v \frac{\partial \rho}{\partial y} -\frac{1}{\rho} \frac{\partial u}{\partial x} -\frac{1}{\rho} \frac{\partial v}{\partial %y}=0
%\]
%and finally multiplying by $\rho^2$ we recover the continuity equation in Eulerian coordinates
\begin{equation} \label{eq:2DEulConsLaw}
	\frac{\partial \rho}{\partial t} + \frac{\partial (\rho u)}{\partial x} + \frac{\partial (\rho v)}{\partial y} =0.
\end{equation}

\noindent
For the momentum equation, we discuss the case $\partial_t w^\Lag = 0$ with the other equation being similar. 
%Similar computations apply for $\partial_t \sigma^\Lag = 0$. Since $w^\Lag=u^\Lag+P_1(\rho^\Lag)$ we have
From 
$$
	\partial_t u^\Lag + P_1^\prime(\rho^\Lag) \partial_t \rho^\Lag = 0,
$$
we compute
%$$
%	\rho \frac{\partial u}{\partial t} + \rho u \frac{\partial u}{\partial x} + \rho v \frac{\partial u}{\partial y} + \rho %P_1^\prime(\rho) \frac{\partial \rho}{\partial t} + \rho u P_1^\prime(\rho) \frac{\partial \rho}{\partial x} + \rho v %P_1^\prime(\rho) \frac{\partial \rho}{\partial y} = 0.
%$$
%Using $\partial_t \rho = - \partial_x (\rho u) - \partial_y (\rho v)$ we obtain
$$
	\rho \frac{\partial u}{\partial t} + \rho u \frac{\partial u}{\partial x} + \rho v \frac{\partial u}{\partial y} - \rho^2 P_1^\prime(\rho) \frac{\partial u}{\partial x} - \rho^2 P_1^\prime(\rho) \frac{\partial v}{\partial y} = 0.
$$
%Now, notice that, adding and subtracting $u\partial_t\rho$, 
%$$
%	\rho \frac{\partial u}{\partial t} + \rho u \frac{\partial u}{\partial x} = \frac{\partial \rho u}{\partial t} + %\frac{\partial \rho u^2}{\partial x} + u \frac{\partial \rho v}{\partial y}.
%$$
and finally obtain the following conservative form: 
%get
%$$
%	\frac{\partial \rho u}{\partial t} + \frac{\partial \rho u^2}{\partial x} + \frac{\partial \rho u v}{\partial y} - %\rho^2 P_1^\prime(\rho) \frac{\partial u}{\partial x} - \rho^2 P_1^\prime(\rho) \frac{\partial v}{\partial y} = 0
%$$
%which can be written in conservative form as
\begin{equation} \label{eq:2DEulSpeedX}
	\frac{\partial \rho w}{\partial t} + \frac{\partial \rho u w}{\partial x} + \frac{\partial \rho v w}{\partial y} = 0.
\end{equation}

%Similarly, for $\partial_t \sigma^\Lag=0$ we get
%\begin{equation} \label{eq:2DEulSpeedY}
%	\frac{\partial \rho \sigma}{\partial t} + \frac{\partial \rho u \sigma}{\partial x} + \frac{\partial \rho v \sigma}%{\partial y} = 0.
%\end{equation}

\noindent
Finally, the system of macroscopic equations in Eulerian coordinates writes as follows:
\begin{equation} \label{eq:2DARZEul}
\begin{aligned}
	\partial_t \rho + \partial_x (\rho u) + \partial_y (\rho v) &=0, \\
	\partial_t (\rho w) + \partial_x (\rho u w) + \partial_y (\rho v w) &= 0, \\
	\partial_t (\rho \sigma) + \partial_x (\rho u \sigma) + \partial_y (\rho v \sigma) &= 0.
\end{aligned}
\end{equation}

\noindent
We end the section with the following remarks on the 2D ARZ-type model~\eqref{eq:2DARZEul}.

\begin{remark}
	If $v=\vref=0$, we recover the one-dimensional ARZ model~\eqref{eq:1D:ARZ}.
\medskip
\noindent
	If we substitute the continuity equation in the two equations of the speeds, we get
	\begin{gather*}
		\partial_t w + u \partial_x w + v \partial_y w = 0,\\
		\partial_t \sigma + u \partial_x \sigma + v \partial_y \sigma = 0,
	\end{gather*}
	which are transport equations for $w$ and $\sigma$.

\medskip 
\noindent
	Recall that $w=u+P_1(\rho)$ and $\sigma=v+P_2(\rho)$. Let $w=c_1$ and $\sigma=c_2$ be constants. Then the system of three equations reduces to
	$$
		\partial_t \rho + \partial_x \left(\rho(c_1-P_1(\rho))\right) +  \partial_y \left(\rho(c_2-P_2(\rho))\right) = 0.
	$$
	Assuming that $\gamma_1=1$ and taking $c_1=\uref$, we obtain the Greenshields' law for the $x$ direction:
	$$
		u = \uref (1-\rho).
	$$
	Changing the values of $c_1$ or $\gamma_1$ will result in several diagrams, as already studied in~\cite{seibold2013NHM}.

\medskip
\noindent
On the other hand, for the $y$ direction we have
	$$
		v = c_2 - P_2(\rho) = c_2 - \frac{\vref}{\gamma_2} \rho^{\gamma_2}, 
	$$
	where $\vref<0$ is the reference velocity in $y$, and it is negative as suggested by the experimental data. Furthermore, in this case, if $c_2=\vref$ and $\gamma_2=1$ we get a linear function such that $v=\vref$ for $\rho=0$ and $v=0$ for $\rho=1$. Actually, as in~\cite{HmVg2017}, the parameters $c_2$ and $\gamma_2$ can be chosen to fit the data and they yield different diagrams, shifted from the naive one if we modify $c_2$, or with a different shape if we modify $\gamma_2$.

\medskip
	The above considerations mean that model~\eqref{eq:2DARZEul} is able to feature different speed diagrams, depending on the modeling of the pressure functions $P_1$ and $P_2$. However, it is worth stressing that this modeling is still done at the macroscopic level, i.e. by means of a given heuristic closure on $P_1$ and $P_2$ as functions of the density. In other words, the different speed diagrams are not derived from microscopic interactions between vehicles, like for kinetic models of traffic flow. For a detailed discussion on these aspects, we refer the reader to recent studies on speed diagrams in one-dimensional~\cite{FermoTosin14,hertyillner12,IllnerKlarMaterne,KlarWegener96,PgSmTaVg,PgSmTaVg3,PgSmTaVg2,VgHmPgTa} and two-dimensional~\cite{HmTaVgZm} kinetic models for traffic flow.
\end{remark}

\section{Properties of the two-dimensional ARZ model} \label{sec:2DMacroModel.properties}
For simplicity of the following discussion, we consider the following non-conservative form of the 2D Aw-Rascle model:

\begin{gather}
	\partial_t \rho + \partial_x (\rho u) + \partial_y (\rho v) =0, \label{s1eq1} \\
	(\partial_t + u\partial_x  + v\partial_y)(u + P_1(\rho)= 0, \label{s1eq2}\\
	(\partial_t + u\partial_x  + v\partial_y)(v + P_2(\rho))= 0. \label{s1eq3}
\end{gather}

\noindent
Using the continuity equation (\ref{s1eq1}) in equations (\ref{s1eq2})-(\ref{s1eq3}), the above system writes, in matrix notations with $U=(\rho,u,v)^T$, as follows

%\begin{gather}
%\partial_t \rho + \rho \partial_x u + u \partial_x \rho + \rho \partial_y v + v \partial_y \rho =0,\label{s2eq1}  \\
%\partial_t u + (u-\rho P^\prime_1(\rho)) \partial_x u + v \partial_y u - \rho P_1^\prime(\rho) \partial_y v = 0, %\label{s2eq2} \\
%\partial_t v + u \partial_x v -\rho P^\prime_2(\rho) \partial_x u + (v - \rho P_2^\prime(\rho)) \partial_y v = 0. %\label{s2eq3} 
%\end{gather}
%using matrix notations, the system (\ref{s2eq1})-(\ref{s2eq3}) writes

\begin{equation}\label{Mateq1}
	\partial_t U + A(U) \partial_x U + B(U) \partial_y U = 0,
\end{equation}
where 
\[ A(U) = \begin{pmatrix}
    	u & \rho & 0\\
        0 & u-\rho P_1^\prime(\rho) & 0\\
        0 & v & -\rho P_1^\prime(\rho)
    \end{pmatrix}
\; \; \text{and} \; \;
	B(U) = \begin{pmatrix}
    	v & 0 & \rho\\
        0 & -\rho P_2^\prime(\rho) & u\\
        0 & 0 & v-\rho P_2^\prime(\rho)
    \end{pmatrix}
\]
\noindent
A system of the form (\ref{Mateq1}) is said to be hyperbolic if for any $\xi=(\xi_1, \xi_2) \in \mathbb{R}^2$, the matrix $C(U, \xi) = \xi_1 A + \xi_2 B$ is diagonalizable \cite{LawrenceBook}. This is the case here and the eigenvalues of the matrix $C(U, \xi)$ are 
\[
\lambda_1 = - (\xi_1 \rho P_1'(\rho) + \xi_2 \rho P_2'(\rho)) \leq \lambda_2 =\xi_1(u - \rho P_1'(\rho)) + \xi_2(v - \rho P_2'(\rho)) \leq \lambda_3 = \xi_1 u + \xi_2 v.
\]
The eigenvalues of $C(U, \xi)$ for $|\xi| = 1$ correspond to the wave speeds and the their corresponding eigenvectors are respectively
\[
r_1 = \begin{pmatrix}
    	\frac{-(u+v)}{v(P_1'(\rho) + P_2'(\rho))}\\
        u/v \\
        1 
    \end{pmatrix}, \; \; \; \;
    r_2 = \begin{pmatrix}
    	0\\
        -1 \\
        1
    \end{pmatrix}, \; \; \; \text{and} \; \; \;
    r_3 = \begin{pmatrix}
    	1\\
        0 \\
        0 
    \end{pmatrix}.
\]

\begin{remark}
	Observe that, as in the classical Aw-Rascle-Zhang model~\cite{aw2000SIAP,ZhangMacro}, the system is not hyperbolic in the presence of a vacuum state, i.e. if $\rho=0$.
    
	\medskip
    
	Moreover, the characteristic wave speeds of system~\eqref{Mateq1} are all less than or equal to the speed $\sqrt{u^2+v^2}$ of vehicles. This supports the anisotropic property of the 2D model, which is known to be a fundamental characteristic for the consistency of second order macroscopic models for traffic flow, cf.~\cite{aw2000SIAP,Daganzo}.
\end{remark}

The second and third eigenvalues, $\lambda_2$ and $\lambda_3$, are linearly degenerate since $\nabla \lambda_2 \cdot r_1 = \nabla \lambda_3 \cdot r_2 =0$ (here $\nabla := (\partial_\rho, \partial_u, \partial_v)$).
For the first eigenvalue, we have 
\[\nabla \lambda_1 \cdot r_1 = \frac{[P_1'(\rho) + P_2'(\rho) + \rho(P_1''(\rho) + P_2''(\rho))](u + v)}{v(P_1'(\rho) + P_2'(\rho) )}.
\]

Clearly, there exist $(\rho, u, v) \in \mathbb{R}^*_+ \times \mathbb{R}^*_+ \times \mathbb{R}^*_+$ such that $\nabla \lambda_1 \cdot r_1 \neq 0$. Therefore $\lambda_1$ is not linearly degenerate. For the other fields to be genuinely nonlinear we need to show that $\nabla \lambda_1 \cdot r_1 \neq 0$ for all $(\rho, u, v)\neq (0, 0, 0)$. This condition is fulfilled provided that the functions $P_1(\rho)$ and $P_2(\rho)$ are such that:
\[
\partial_\rho (\rho P_1'(\rho)) \neq \partial_\rho (\rho P_2'(\rho)). 
\]
The above condition is satisfied as long as the reference velocities $U_{ref}$ and $V_{ref}$, introduced in (\ref{functionsP}), are different. Since $\lambda_1$ is genuinely nonlinear, for $U_{ref} \neq V_{ref}$, then the associated waves are either shocks or rarefaction waves. 

The Riemann invariants, in the sense of Lax \cite{LaxBook}, associated with the eigenvalues $\lambda_1, \lambda_2$ and $\lambda_3$ are respectively:
\[ z_1 = u + v + P_1(\rho) + P_2(\rho), \; \; \; z_2 = u + v, \; \; z_3 = 
u.\]

\subsection{An overview of the Riemann problem associated with the system}

%\todo{ It is not clear what the setup for the Riemann problem is. Normally in 2d you would have four sets
%of initial data. However, the discussion below is related to only two sets of initial data. Here, I suggest to 
%state the Riemann problem first, so that it is clear that we deal with in fact 1D %Riemann problems} 

The main idea behind the approximation of the solution to the two-dimensional Riemann problem was the coupling of two-half Riemann problems in one dimension and then use the transition from one lane to another, subject to some admissibility conditions. In this section, we provide a detailed discussion, of the main scenarios under consideration in this study, as well as the suitable Riemann solvers at the interfaces of the lanes in order to capture the two-dimensional aspect of the problem. Following this approach, the solution to the system (\ref{s1eq1})-(\ref{s1eq3}) consists of either a wave of the first family (1-shock or 1-rarefaction) or a wave of the second
family (2-contact discontinuity). The properties of these two families of waves can be summarized as follows.

\paragraph{First characteristic field:} A wave of the first family is generated when a state on left, denoted $U_l=(\rho_l, u_l, v_l)$, is connect to a state on the right, denotes $U_r=(\rho_r, u_r, v_r)$,  through the same Riemann invariant curve associated with the eigenvalue $\lambda_1$, i.e. $z_1(U_l) = z_1(U_r)$. Therefore, we have
\begin{equation}
u_l + v_l + P_1(\rho_l) + P_2(\rho_l) = u_r + v_r + P_1(\rho_r) + P_2(\rho_r).
\end{equation}
We can distinguish the following scenarios:
\begin{itemize}
\item if $u_l + v_l > u_r + v_r$, then this wave (of the first family) is a 1-shock i.e. a jump discontinuity, traveling with the speed
\begin{equation}
s = \frac{\rho_r(u_r+v_r) - \rho_l(u_l + v_l)}{\rho_r - \rho_l}.
\end{equation}

\item if $u_l + v_l < u_r + v_r$, then this wave (of the first family) is a 1-rarefaction, i.e. a continuous solution of the form $(\rho, u, v)(\xi)$ (with $\xi = \frac{f(x, y)}{t}$), given by 

\begin{equation}
\begin{pmatrix}
    	\rho'(\xi)\\
        u'(\xi) \\
        v'(\xi) 
    \end{pmatrix} 
    = \frac{r_1(U(\xi))}{\nabla \lambda_1(U(\xi)) \cdot r_1(U(\xi))}, \; \; \; \text{if} \; \;  \lambda_1(U_l)\leq \xi \leq \lambda_1(U_r),   
\end{equation}
whereas 
\[
(\rho, u, v)(\xi) = 
\begin{cases}
(\rho_l, u_l, v_l) \; \; \; \text{for} \; \; \; \xi <\lambda_1(U_l),\\
(\rho_r, u_r, v_r) \; \;  \text{for} \; \; \; \xi >\lambda_1(U_r).\\
\end{cases}
\]
\end{itemize}

\paragraph{Second characteristic field:}
A wave of the second family i.e. a 2-contact discontinuity is generated when 
$$z_2(U_l) = z_2(U_r) \Longrightarrow u_l + v_l = u_r +v_r.$$
In this case, the contact discontinuity between a state on the left, $U_l=(\rho_l, u_l, v_l)$, and a state on the right, $U_r=(\rho_r, u_r, v_r)$, travels with a speed $\eta = u_l + v_l = u_r +v_r$.

\paragraph{Third characteristic field:}
A wave of the third family i.e. a 3-contact discontinuity is generated when 
$$z_3(U_l) = z_3(U_r) \Longrightarrow u_l  = u_r .$$
In this case, the contact discontinuity between a state on the left, $U_l=(\rho_l, u_l, v_l)$, and a state on the right, $U_r=(\rho_r, u_r, v_r)$, travels with a speed $\eta = u_l  = u_r $.

\subsection{Solution to the Riemann problem associated with the system (\ref{s1eq1})-(\ref{s1eq3})}

This section describes the solutions to the Riemann problem for the system (\ref{s1eq1})-(\ref{s1eq3}), by combining the previously described elementary waves.  Let  $U_l = (\rho_l, u_l, v_l)$ and $U_r= (\rho_r, u_r, v_r)$ be the initial data on the left and on the right, respectively. The solutions to the system (\ref{s1eq1})-(\ref{s1eq3}) for these initial data consist of the following cases.

\paragraph{Case 1.} If $\sigma_l = v_l + P_2(\rho_l) = \sigma_r = v_r + P_2(\rho_r) = 0$, then the model bowls down to the one-dimensional Aw-Rascle model, i.e. the initial data reduced to $U_l=(\rho_l, u_l)$ and $U_r= (\rho_r, u_r)$. Therefore, the state on the left $U_l = (\rho_l, u_l)$ is connect to an intermediate state $U^\star = (\rho^\star, u^\star)$ through a wave of the first family, i.e. either a 1-shock wave if $u_l > u_r$ (which corresponds to braking) or a 1-rarefaction wave if $u_l < u_r$ (which corresponds to an acceleration). Then, the intermediate state $U^\star$ is connected to the state on the right $U_r = (\rho_r, u_r)$ through a wave of the second family, i.e. a 2-contact discontinuity. Therefore, an admissible intermediate state, $U^\star = (\rho^\star, u^\star)$, is defined such that: 

$$u^\star = u_r \; \; \;  \text{and} \; \; \;  \rho^\star = P_1^{-1}( u_l + u^\star + P_1(\rho_l)).$$ 

\paragraph{Case 2.} If $\sigma_l = v_l + P_2(\rho_l) = \sigma_r = v_r + P_2(\rho_r) > 0$, then we have the following admissible scenarios:

\begin{itemize}
\item If $v_l = \sigma_l - P_2(\rho_l) > 0$ (lane changing possibility), then the state on the left $U_l = (\rho_l, u_l, v_l)$ is connected to an intermediate state on the left, $U_l^\star = (\rho_l^\star, u_l^\star, v_l^\star)$,  through a wave of the first family, i.e. either a 1-shock (if $u_l > u_r$) or a rarefaction wave (if $u_l < u_r$). Then, the intermediate state $U_l^\star$ is connected to a vacuum  state through a wave of the second family, i.e. a 2-contact discontinuity.  Therefore, an admissible intermediate state, $U_l^\star =  (\rho_l^\star, u_l^\star, v_l^\star)$, is defined such that: 
$$
\begin{cases}
u_l + P_1(\rho_l) = u_l^\star + P_1(\rho_l^\star), \\
v_l + P_2(\rho_l) = v_l^\star + P_2(\rho_l^\star), \\
u_l^\star = u_l + v_l.
\end{cases}
$$
The vacuum state, is then connected to an intermediate state on the right, $U_r^\star = (\rho_r^\star, u_r^\star, v_r^\star)$, through a wave of the first family, i.e. either a 1-shock (if $u_l^\star > u_r^\star$) or a 1-rarefaction wave (if $u_l^\star < u_r^\star$).  Then, the intermediate state, $U_r^\star$, is connected to the state on the right, $U_r = (\rho_r, u_r, v_r)$, through a wave of the third family, i.e. a 3-contact discontinuity. Therefore, an admissible intermediate state, $U_r^\star = (\rho_r^\star, u_r^\star,  v_r^\star)$, is defined such that:
$$
\begin{cases}
u_l^\star + P_1(\rho_l^\star) = u_r^\star + P_1(\rho_r^\star), \\
v_l^\star + P_2(\rho_l^\star) = v_r^\star + P_2(\rho_r^\star), \\
u_r^\star = u_r.
\end{cases}
$$
 
\item  If $v_l = \sigma_l - P_2(\rho_l) = 0$ (no lane changing possibility), then the state on the left $U_l = (\rho_l, u_l, v_l)$ is connected to an intermediate state on the left, $U_l^\star = (\rho_l^\star, u_l^\star, v_l^\star)$, through a 1-rarefaction wave. Then, the intermediate state $U_l^\star$ is connected to a vacuum  state through a wave of the second family, i.e. a 2-contact discontinuity, and an admissible intermediate state, $U_l^\star =  (\rho_l^\star, u_l^\star, v_l^\star)$, is defined such that: 
$$
\begin{cases}
u_l + P_1(\rho_l) = u_l^\star + P_1(\rho_l^\star), \\
P_2(\rho_l) = v_l^\star + P_2(\rho_l^\star), \\
u_l^\star = u_l.
\end{cases}
$$
The vacuum state, is then connected to an intermediate state on the right, $U_r^\star = (\rho_r^\star, u_r^\star, v_r^\star)$, through a either a 1-shock or a 1-rarefaction wave; whereas the intermediate state, $U_r^\star$, is connected to the state on the right, $U_r = (\rho_r, u_r, v_r)$, through a 3-contact discontinuity. An admissible intermediate state, $U_r^\star = (\rho_r^\star, u_r^\star,  v_r^\star)$, is therefore defined such that:

$$
\begin{cases}
u_l^\star + P_1(\rho_l^\star) = u_r^\star + P_1(\rho_r^\star), \\
v_l^\star + P_2(\rho_l^\star) = v_r^\star + P_2(\rho_r^\star), \\
u_r^\star = u_r.
\end{cases}
$$
\end{itemize}

\begin{figure*}[t!]
    \centering
    \begin{subfigure}[t]{0.5\textwidth}
       \centering
        \includegraphics[width=\textwidth]{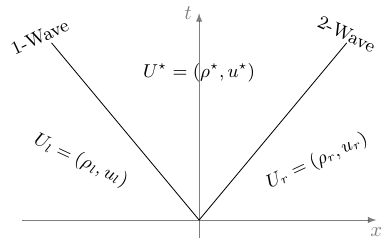}
        \caption{Case 1.}
    \end{subfigure}
    ~
    \begin{subfigure}[t]{0.5\textwidth}
        \centering
        \includegraphics[width=\textwidth]{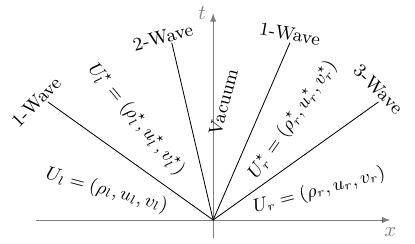}
        \caption{Case 2.}
    \end{subfigure}
    \caption{Illustration of the solutions to the Riemann problem associated with the system (\ref{s1eq1})-(\ref{s1eq3}).}
\end{figure*}

\section{Numerical simulations}\label{sec:numericalSimulations}

In Section \ref{sec:2DARmodel} we have derived the two-dimensional ARZ-type model \eqref{eq:2DARZEul} starting from the two-dimensional FTL-type model \eqref{eq:2DFTL}. This shows the connection between the microscopic model and a semi-discretization of the macroscopic continuum model in Lagrangian coordinates.

Here, we aim to provide also some numerical evidence of the connection between the two models. To achieve this, we will present some simulations obtained using a first-order in time approximation of the particle model \eqref{eq:2DFTL} and a first-order finite volume scheme of the Eulerian version of the macroscopic model \eqref{eq:2DARZEul}. Furthermore, the numerical results presented in this section used a zero-flux boundary in $y$ direction. However, a reflecting boundary conditions might be an alternative option.

\subsection{Description of the schemes}

\paragraph{Microscopic model:} For the microscopic model \eqref{eq:2DFTL}, we use the explicit Euler scheme in time, so that we obtain the following fully-discretized equations for the update of the positions and the speeds at time $t^{n+1}$
\begin{align*}
	x_i^{n+1} &= x_i^n + \Delta t \, u_i^n,\\
	y_i^{n+1} &= y_i^n + \Delta t \, v_i^n,\\
	u_i^{n+1} & = u_i^n + \Delta t \, C_1 \left( \frac{u_{j(i)}^n-u_i^n}{\left(x_{j(i)}^n-x_i^n\right)\Delta A^{\gamma_1,n}}  + \frac{v_{j(i)}^n-v_i^n}{\left(y_{j(i)}^n-y_i^n\right)\Delta A^{\gamma_1,n}}\right),	\\
	v_i^{n+1} & = v_i^n + \Delta t \, C_2 \left( \frac{u_{j(i)}^n-u_i^n}{\left(x_{j(i)}^n-x_i^n\right)\Delta A^{\gamma_2,n}}  + \frac{v_{j(i)}^n-v_i^n}{\left(y_{j(i)}^n-y_i^n\right)\Delta A^{\gamma_2,n}}\right),
\end{align*}
where $\Delta t$ is the time-step. The choice of the initial conditions as well as the computation of the characteristics of the interacting vehicles will be discussed later.

\paragraph{Macroscopic model:} For the macroscopic model in Eulerian coordinates \eqref{eq:2DARZEul}, we consider a classical first-order finite volume approximation. We divide the domain $(x,y) \in [a^x,b^x]\times [a^y,b^y]$ into $N^x\times N^y$ cells
\[
	I_{ij} = [x_{i-1/2},x_{i+1/2}] \times [y_{j-1/2},y_{j+1/2}], \quad i=1,\dots,N^x, \ j=1,\dots,N^y,
\]
where $x_{i+1/2}-x_{i-1/2}=\Delta x$, $y_{j+1/2}-y_{j-1/2}=\Delta y$ and the mid-points are	$x_i$, $i=1,\dots,N^x$, and $y_j$, $j=1,\dots,N^y$. The length and the width of the road are $L^x: = b^x-a^x$ and $L^y: = b^y-a^y$, respectively. Let us to denote by
\begin{gather*}
	\overline{\rho}_{ij}(t) = \frac{1}{\Delta x \Delta y} \iint_{I_{ij}} \rho(t,x,y) \mathrm{d}x \mathrm{d}y, \quad
	\overline{\rho w}_{ij}(t) = \frac{1}{\Delta x \Delta y} \iint_{I_{ij}} (\rho w)(t,x,y) \mathrm{d}x \mathrm{d}y, \\
	\overline{\rho \sigma}_{ij}(t) = \frac{1}{\Delta x \Delta y} \iint_{I_{ij}} (\rho \sigma)(t,x,y) \mathrm{d}x \mathrm{d}y
\end{gather*}
the cell averages of the exact solution of the system \eqref{eq:2DARZEul} at time $t$. Moreover, in order to abbreviate the notation, we define $q:=(\rho,\rho w,\rho\sigma)$ and we denote by $\overline{Q}_{ij}(t)$ the vector of the numerical approximation of the cell-averages $\overline{q}_{ij}(t):=(\overline{\rho}_{ij}(t),\overline{\rho w}_{ij}(t),\overline{\rho\sigma}_{ij}(t))$. Then, system \eqref{eq:2DARZEul} can be then rewritten as
\begin{equation} \label{eq:2Dvectorial}
	\partial_t q(t,x,y) + \partial_x f(q(t,x,y)) + \partial_y g(q(t,x,y)) = 0,
\end{equation}
where
\[
	f(q) = \begin{pmatrix}
	\rho u\\
	\rho u w\\
	\rho u \sigma
	\end{pmatrix},
	\quad
	g(q) = \begin{pmatrix}
	\rho v\\
	\rho v w\\
	\rho v \sigma
	\end{pmatrix}.
\]

By integrating equation \eqref{eq:2Dvectorial} over a generic cell $I_{ij}$ of the grid, dividing by $\Delta x \Delta y$ and finally using the explicit Euler scheme with time-step $\Delta t$, we get the fully-discrete scheme for the approximation of the solution at time $t^{n+1}=t^n+\Delta t$:
\begin{equation} \label{eq:2DARZscheme}
\overline{Q}_{ij}^{n+1} = \overline{Q}_{ij}^n - \frac{\Delta t}{\Delta x} \left[ F_{i+1/2,j} - F_{i-1/2,j} \right] - \frac{\Delta t}{\Delta y} \left[ G_{i,j+1/2} - G_{i,j-1/2} \right].
\end{equation}
where
\begin{gather*}
	F_{i+1/2,j} = \mathcal{F}\left(q(t^n,x_{i+1/2}^-,y_j),q(t^n,x_{i+1/2}^+,y_j)\right) \approx f(q(t^n,x_{i+1/2},y_j))\\
	G_{i,j+1/2} = \mathcal{G}\left(q(t^n,x_i,y_{j+1/2}^-),q(t^n,x_i,y_{j+1/2}^+)\right) \approx g(q(t^n,x_i,y_{j+1/2}))
\end{gather*} 
are the numerical fluxes defined by $\mathcal{F}$ and $\mathcal{G}$, being approximate Riemann solvers (in the following the local Lax-Friedrichs). We observe that in order to compute the numerical flux, we need to to know the solution in the mid points of the four boundaries of a cell $I_{ij}$. If we use a first-order scheme, then
\begin{gather*}
	q(t,x_{i+1/2}^-,y_j)\approx\overline{Q}_{ij}^n, \quad q(t,x_{i+1/2}^+,y_j)\approx\overline{Q}_{i+1,j}^n,\\
	q(t,x_i,y_{j+1/2}^-)\approx\overline{Q}_{ij}^n, \quad q(t,x_i,y_{j+1/2}^+)\approx\overline{Q}_{i,j+1}^n.
\end{gather*} 

In the next section we will show that using scheme \eqref{eq:2DARZscheme} will lead to the same results obtained using the above discretized particle model with a large number of vehicles.

\begin{remark}
	We point out that the discretization introduced above has some close relation to multilane models if one uses a very coarse discretization in $y$-direction, e.g. if the number of discretization points is equal to the number of lanes. However, the model introduced here is based on the idea to treat the lanes as continuum. Therefore, if we use a coarse discretization in $y$-direction, this will introduce a strong numerical diffusion to the solution that might lead to qualitatively different solutions compared to those predicted by a (theoretical) model.
    Moreover, a coarse discretization of the equations lead to a conservative model which contrasts with multi-lane models with exchange terms on the right-hand side.
\end{remark}

\subsection{Examples}

In the subsequent sections, we will present three numerical simulations. In the first one, we compare the microscopic model and the above first-order scheme for the macroscopic model using a basic test problem, which generalizes in two dimensions the problem proposed in \cite{aw2002SIAP} and recalled in Section \ref{sec:preliminary}. Then, we propose two simulations regarding the macroscopic model only, in order to show that it is capable to reproduce typical situations in traffic flow.

In all cases, we consider the following settings:
\begin{gather*}
	[0,1] = 1~\text{km}, \quad L^y = 0.012, \quad \rho_{\max} = 1, \quad U_\text{ref} = 1, \quad V_\text{ref} = 0.009, \quad \gamma_1 = \gamma_2 = 1.
\end{gather*}
We use free-flow boundary conditions in the direction of the flow of vehicles, while we use zero-flux boundary conditions at the edges of the road section so that density cannot neither enter nor leave.

\paragraph{Connection between Micro and Macro models.} Here, we will numerically show that for $\Delta t$, $\Delta X$ and $\Delta Y$ tending to $0$ one obtains an approximation of the two-dimensional system \eqref{eq:2DARZEul} in Eulerian coordinates. From the particle point of view, this means that the number of vehicles should increase in order to get the desired approximation of the macroscopic equations.

From the macroscopic point of view, we consider a Riemann Problem given by the following initial conditions
\[
	\rho(0,x,y)=0.05, \quad u(0,x,y)=\begin{cases} 0.8, & x\geq 0\\0.05 & x<0\end{cases} \quad v(0,x,y)=\begin{cases} -0.001, & y\geq L^y/2\\0.001 & y<L^y/2\end{cases},
\]
which define four states with discontinuity located along $x = 0$ and $y = L^y/2$. In other words, we are assuming that all vehicles in the left part of the road are traveling towards the right part, while vehicles in the right part are traveling towards the left part. We denote the four states by NE (north-east), NW (north-west), SW (south-west) and SE (south-east). Thus, we have
\begin{gather*}
	\rho_\text{NE} = \rho_\text{NW} = \rho_\text{SE} = \rho_\text{SW} = 0.05,\\
	(\rho u)_\text{NE} = (\rho u)_\text{SE} = 0.04, \quad (\rho u)_\text{NW} = (\rho u)_\text{SW} = 0.0025,\\
	(\rho v)_\text{NE} = (\rho v)_\text{NW} = -5\times 10^{-5}, \quad (\rho v)_\text{SE} = (\rho v)_\text{SW} = 5\times 10^{-5}.
\end{gather*}

The discretization size is chosen as $\Delta x = \Delta X = \frac{1}{200}$ and $\Delta y = \Delta Y = \frac{L^y}{32}$. From the particle point of view, since the density is constant and equal to $0.05$, this yields $320$ cars per kilometer. The initial conditions for the particle model are assigned as follows. We focus on the simplest case in which $N^y = 4$ and we use two indices for labeling the microscopic states of vehicles: the first counts the vehicles, while the second takes into account the ``lane''. We put the same number of vehicles in each of the four lanes and then proceed as follows:
\begin{enumerate}
	\item firstly, we choose the initial position of the first vehicle belonging to lane $1$ and lane $3$, so that
	\[
		x_{1,1}(0) = x_{1,3}(0) = a^x,
	\]
	and we compute the position at initial time of the first vehicle in lane $2$  such that the density \eqref{eq:2Ddensity} is equal to $0.05$. Thus
	\[
		x_{1,2}(0) = x_{1,1}+\frac{\Delta X\Delta Y}{0.05\left|y_{1,2}-y_{1,1}\right|} = x_{1,3}+\frac{\Delta X\Delta Y}{0.05\left|y_{1,2}-y_{1,3}\right|}.
	\]
	The initial position of the first vehicle in lane $4$ is finally chosen as $x_{1,4}(0) =x_{1,2}(0)$.
	\item Let $d = x_{1,2} - x_{1,1}$, then all the positions at initial time are given by:
	\[
		x_{i,j} (0) = x_{i-1,j} (0) + 2\,d, \quad i=2,\dots, \ j=1,2,3,4.
	\]
	\item The initial speeds in $x$-direction are:
	\[
		u_{i,j}(0) = \begin{cases}
		0.05 & \text{if $x_{i,j}(0)\leq 0$},\\
		0.8 & \text{if $x_{i,j}(0)>0$}.
		\end{cases}
	\]	
	\item Finally, the initial speed in $y$-direction are:
	\[
		v_{i,j}(0) = \begin{cases}
		0.001 & \text{if $j=1,2$},\\
		-0.001 & \text{if $j=3,4$}.
		\end{cases}
	\]
\end{enumerate}

The above artificial initial conditions for the particle model are induced by the fact that vehicles in the right half part of the road travel with positive lateral speeds, while vehicles in the left half part of the road travel with negative lateral speeds. Thus, it is quite natural to assume that vehicles in the first lane interact with vehicles in the second lane, as well as vehicle in the third lane. While vehicles in the second and the fourth lane interact with vehicles in the third lane.

\begin{figure}
	\centering
	\includegraphics[width=\textwidth]{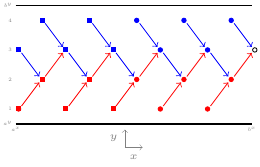}
	\caption{Initial condition for the particle model.\label{fig:MicroInitCond}}
\end{figure}

\medskip

This situation is depicted in Figure \ref{fig:MicroInitCond} for the case of four lanes and it can be easily generalized to the case of an arbitrary number of lanes. Dots and squares represent the vehicles on the road. More precisely, the red cars are in the first and second lane, traveling with a positive lateral speed (i.e., towards the left part of the road). On the other hand, the blue cars are in the third and fourth lane, traveling with a negative lateral speed (i.e., towards the left part of the road). The squares and the dots identify the vehicles having speed $0.05$ and $0.8$ in $x$-direction, respectively. The arrows show, for each car, the interacting vehicle. The empty black circle is the ``ghost'' car, i.e. the boundary condition, which is necessary to compute the microscopic states of the last cars in lane $2$ and $4$. The positions and the speeds of the ghost car are updated at each time using the positions and the speeds of the last car in the same lane, i.e. the third one.

\begin{figure}[t!]
	\centering
	\includegraphics[width=0.49\textwidth]{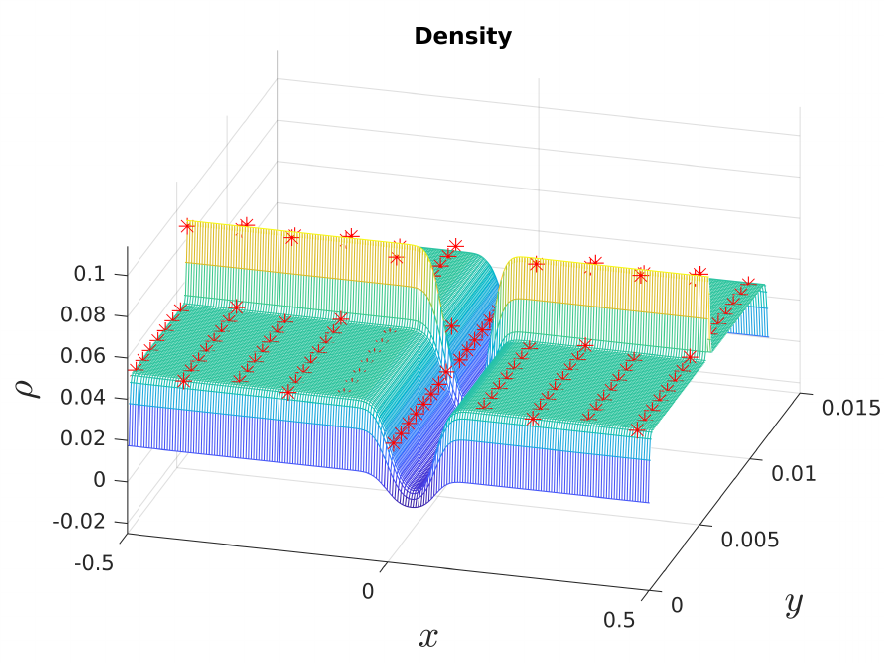}
	\includegraphics[width=0.49\textwidth]{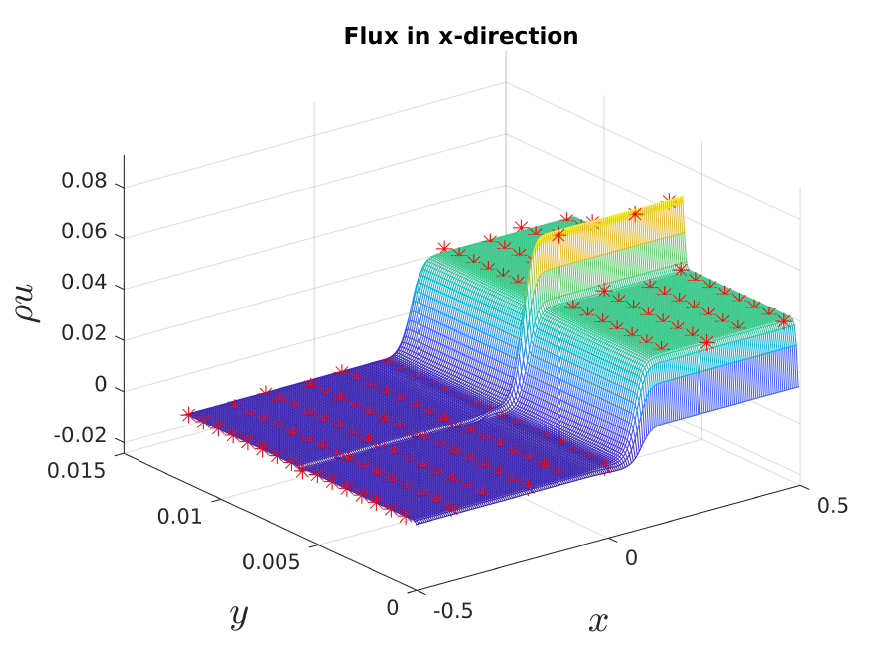}
	\includegraphics[width=0.49\textwidth]{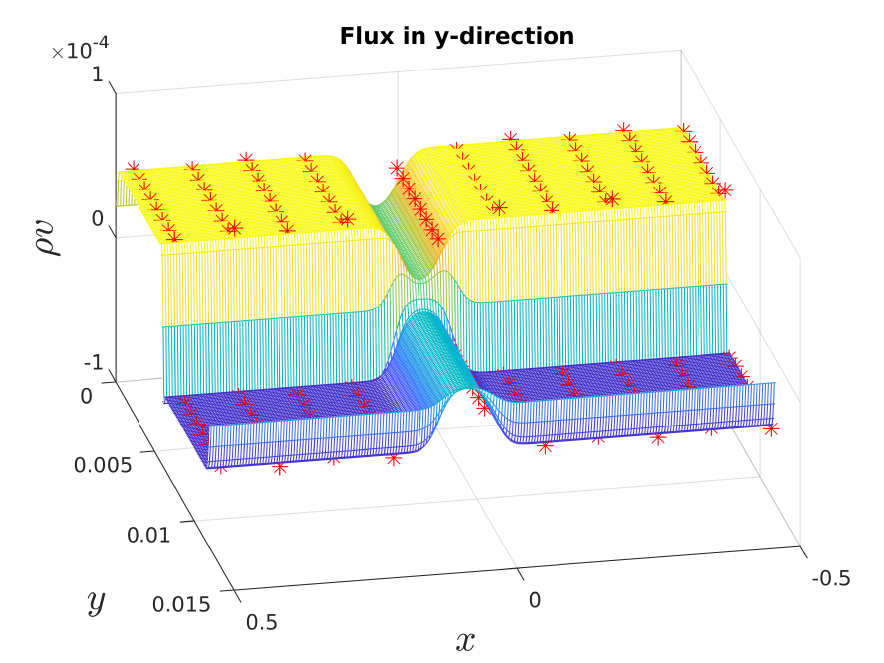}
	\caption{Density $\rho$ (top-left), flux $\rho u$ (top-right) and flux $\rho v$ (bottom-center) profiles at time $T_\text{fin} = 0.1$ provided by the two-dimensional second-order macroscopic model~\eqref{eq:2DARZEul}. The red $*$-symbols shows the values of the density around each car and of the fluxes provided by the two-dimensional microscopic model~\eqref{eq:2DFTL}.\label{fig:SimulationMicroMacro}}
\end{figure}

\medskip

In Figure \ref{fig:SimulationMicroMacro}, we present the density and the fluxes profiles at $T_\text{fin} = 0.1$ provided by the macroscopic model. The time step is chosen according to the CFL condition. We consider a small final time in order to guarantee that, in the particle model, vehicles remain in their part of the domain given at initial time. At each time step, the density around each car is computed again using equation \eqref{eq:2Ddensity} and for each car the interacting vehicle is chosen by \eqref{eq:FieldCar}. The values of the density around each car and of the fluxes provided by the particle model are shown with red $*$-symbols. We notice that the microscopic and the macroscopic model seem to produce the same profiles at final time. Moreover, the density has the same decrease in the center of the road as observed in the 1D simulation, see Figure~\ref{fig:1D}. This happens because the initial condition in $x$ is similar to that given in Section \ref{sec:preliminary}. Finally, we observe that the density tends to $0$ in the left-most and in the right-most part of the road because of the initial condition on the lateral velocity, which assumes that the flow is going towards the center part of the road.

The following two examples concern only the macroscopic model, and they aim to show that a such model reproduces typical traffic flow situations, including overtaking scenarios.

\paragraph{Going to the left.}

First, we propose the case of an overtaking to the left. Thus, using the same label introduced in the previous paragraph, we consider the following initial condition for the density and the speeds:
\begin{gather*}
	\rho_\text{NE} = 0.05, \quad \rho_\text{NW} = \rho_\text{SE} = 0.4, \quad \rho_\text{SW} = 0.6,\\
	u_\text{NE} = u_\text{NW} = 0.8, \quad u_\text{SW} = 0.65, \quad u_\text{SE} = 0.35,\\
	v_\text{NE} = v_\text{NW} = v_\text{SE} = 0, \quad v_\text{SW} = 0.004.
\end{gather*}

In other words, we are taking into account an initial situation in which the flow in the SW region of the road is higher than the flow in the SE region. In a simple one-dimensional model, we expect to have a backward propagating shock since vehicles cannot go in the SE region freely. Using a two-dimensional model, instead, we can observe the overtaking: vehicles in the SW region move in the NE region in order to avoid the slow mass ahead. The initial lateral speed in the SW zone is chosen to be positive in order to speed up the overtaking.

\begin{figure}[t!]
	\centering
	\includegraphics[width=0.32\textwidth]{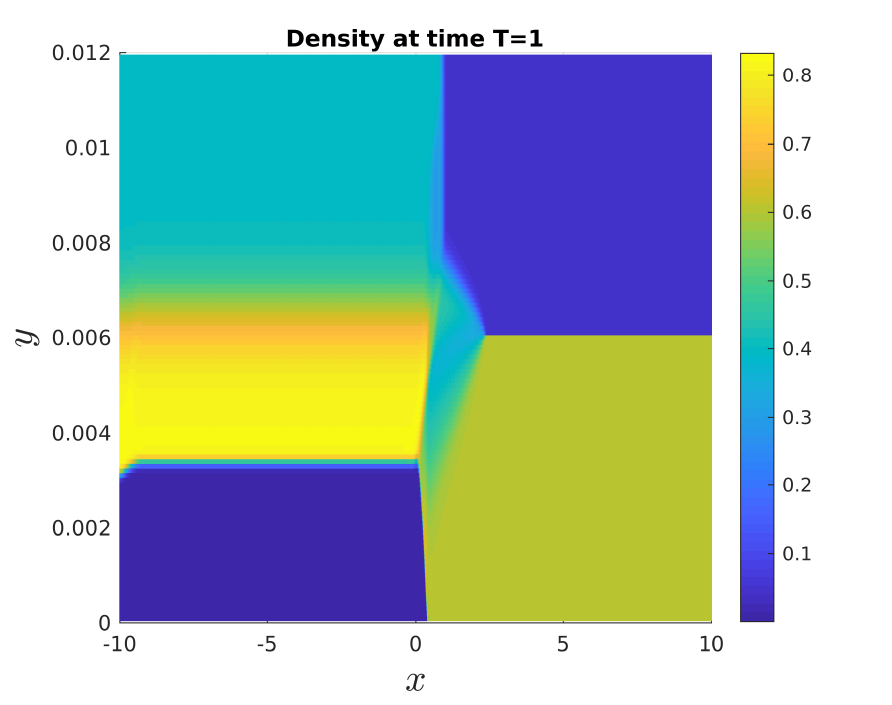}
	\includegraphics[width=0.32\textwidth]{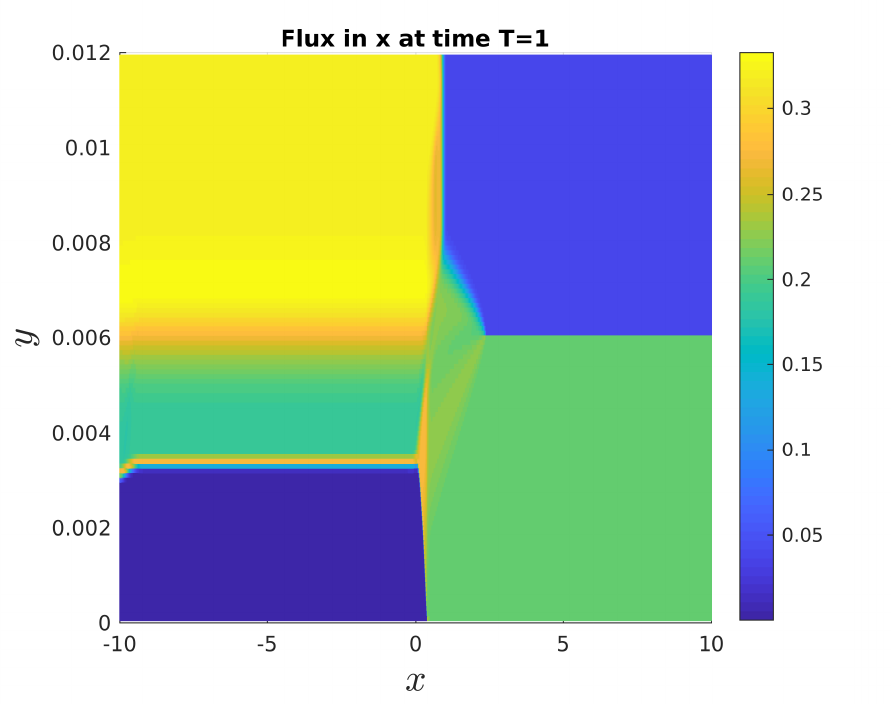}
	\includegraphics[width=0.32\textwidth]{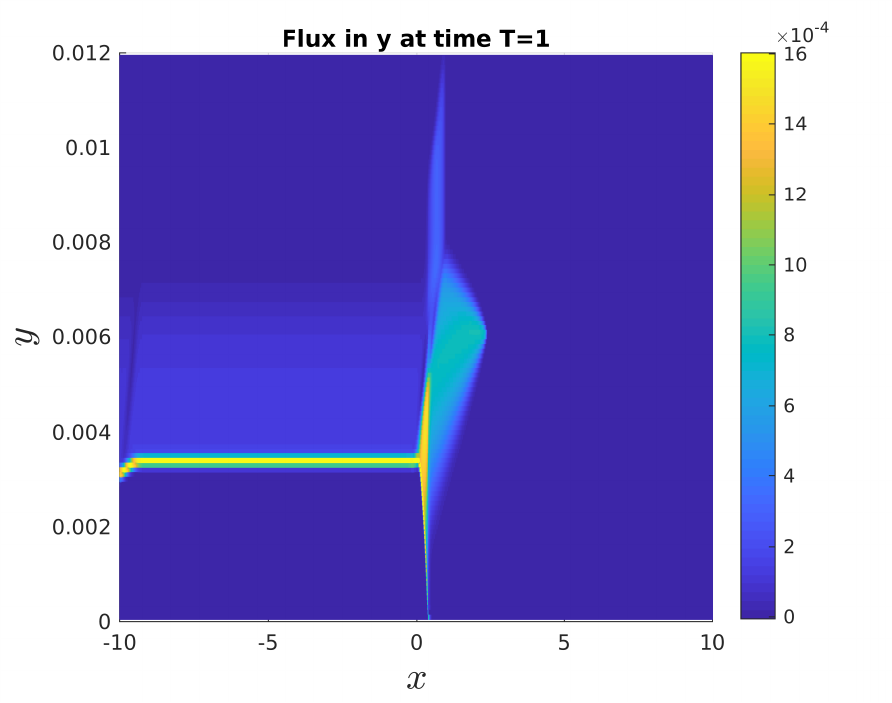}
	\includegraphics[width=0.32\textwidth]{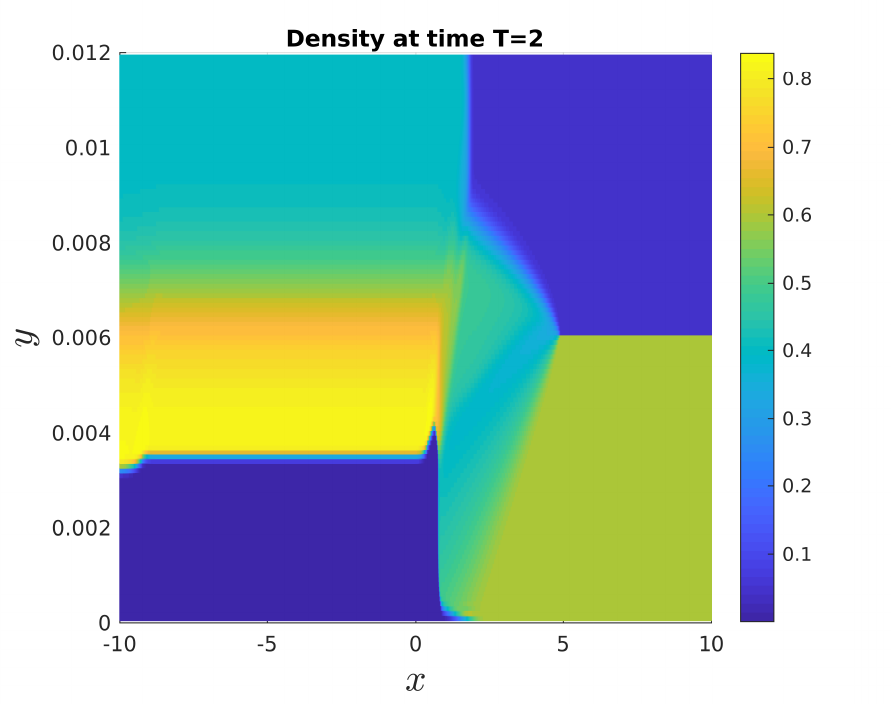}
	\includegraphics[width=0.32\textwidth]{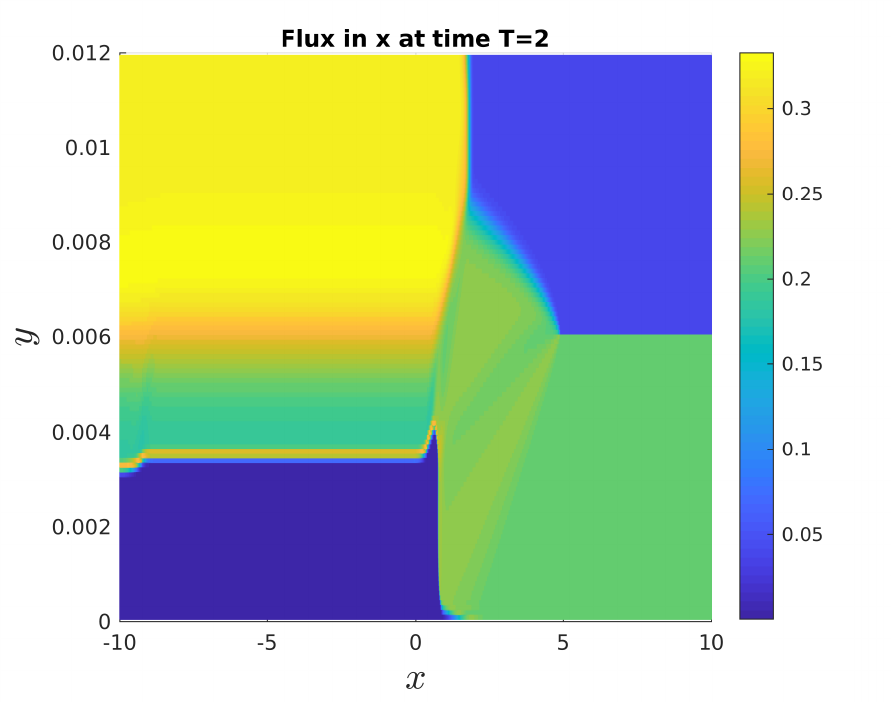}
	\includegraphics[width=0.32\textwidth]{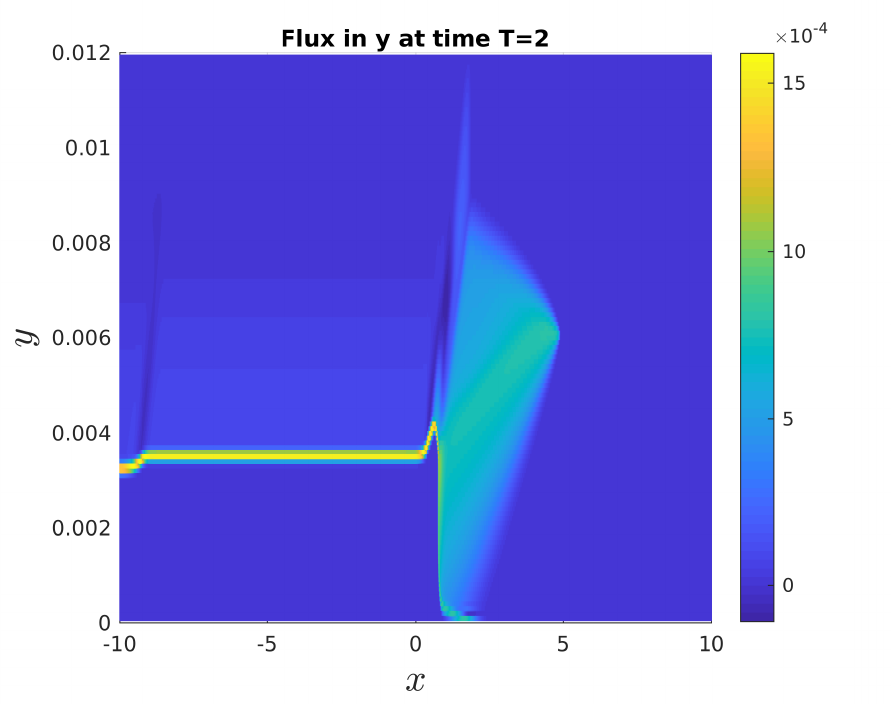}
	\includegraphics[width=0.32\textwidth]{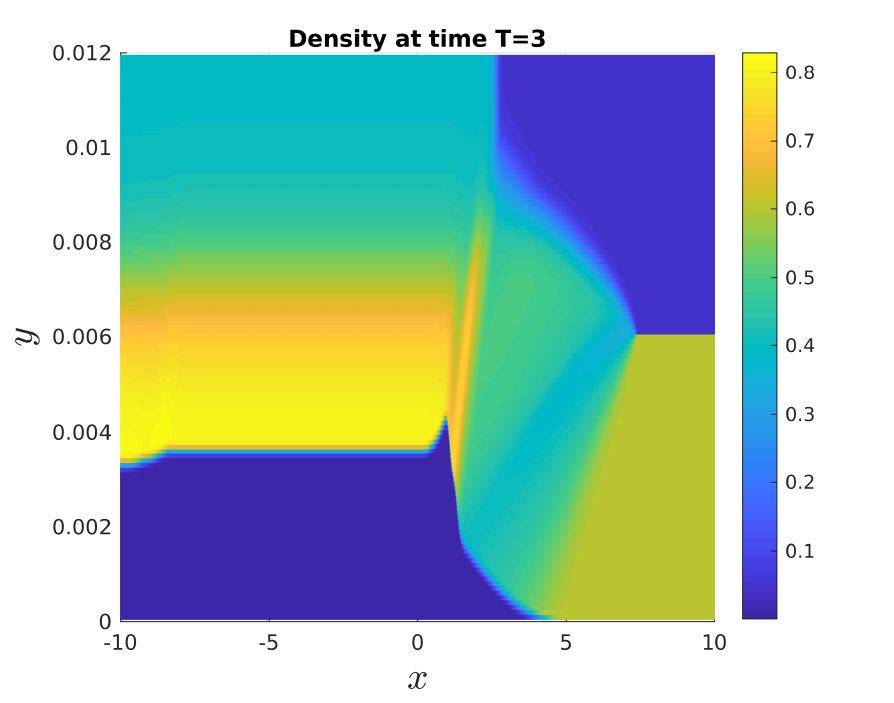}
	\includegraphics[width=0.32\textwidth]{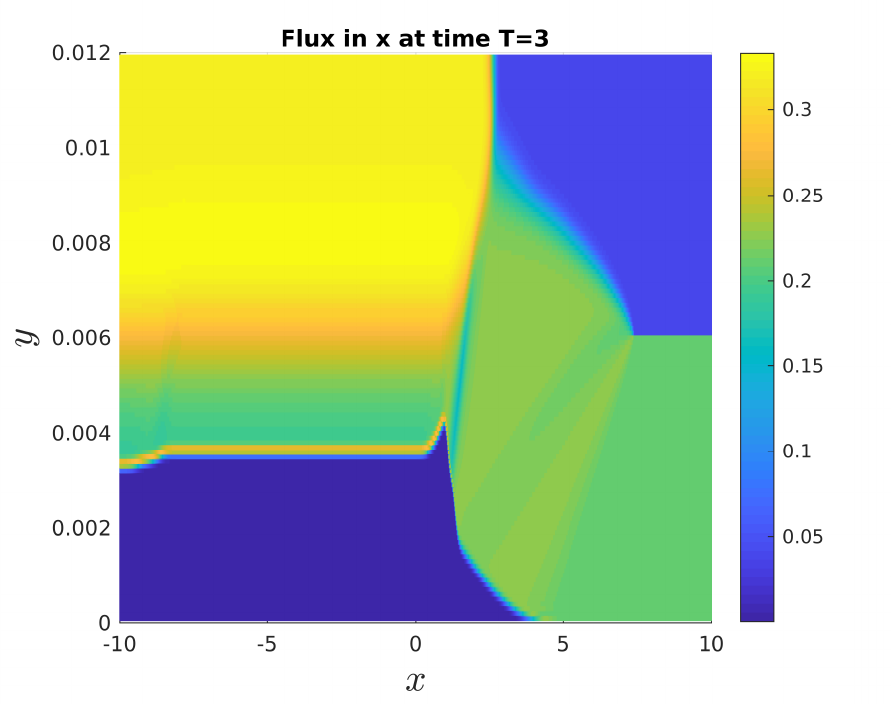}
	\includegraphics[width=0.32\textwidth]{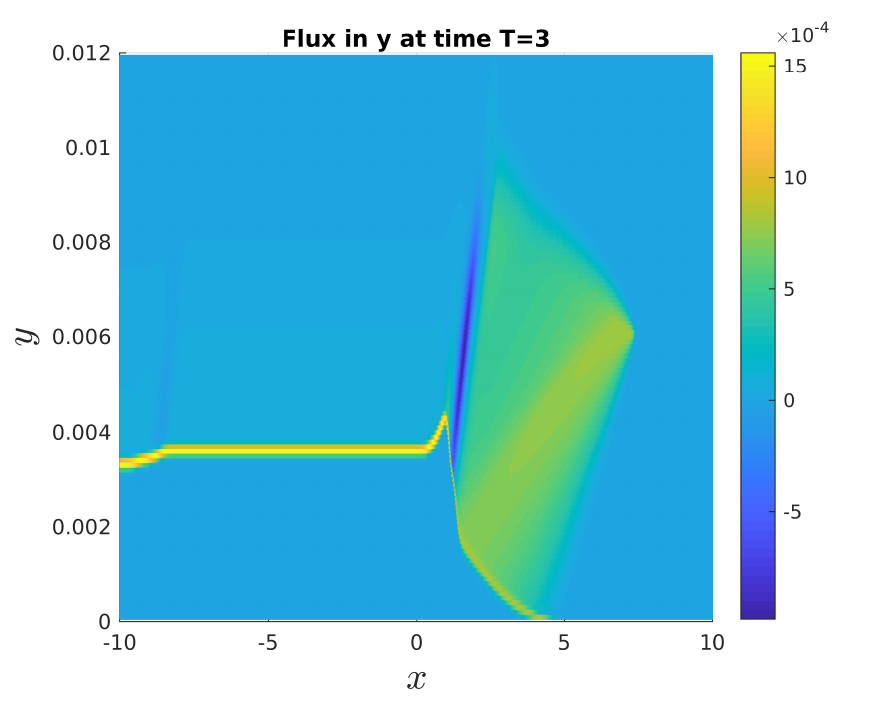}
	\caption{Density $\rho$ (left column), flux $\rho u$ (center column) and flux $\rho v$ (right column) profiles provided by the two-dimensional second-order macroscopic model~\eqref{eq:2DARZEul} at time $T_\text{fin} = 1$ (top row), $T_\text{fin} = 2$ (center row) and $T_\text{fin} = 3$ (bottom row).\label{fig:ToLeft}}
\end{figure}

In Figure \ref{fig:ToLeft} we show the density and the flux profiles provided by the 2D ARZ-type model at time $T=1$ (top row), $T=2$ (center row) and $T=3$ (bottom row).

\paragraph{Going to the right.}

Finally, we study an opposite scenario in which the overtaking is to the right.  We consider the following initial conditions for the density and the speeds:
\begin{gather*}
\rho_\text{NE} = 0.9, \quad \rho_\text{NW} = 0.7, \quad \rho_\text{SE} = \rho_\text{SW} = 0.05,\\
u_\text{NE} = 0.1, \quad u_\text{NW} = 0.7, \quad u_\text{SW} = u_\text{SE} = 1,\\
v_\text{NE} = v_\text{NW} = v_\text{SE} = v_\text{SW} = 0.
\end{gather*}

Now, vehicles in the NW region would travel towards the SE region in order to overtake the slower mass in the NE zone. However, notice that, in contrast to the previous example, the lateral speed is zero everywhere. In fact, our aim is to show that the overtaking takes place also if the initial lateral speed is zero. The macroscopic equations force the lateral speed of vehicle being in the NW region to become negative and thus cars can overtake and travel in the right part of the road.

\begin{figure}[t!]
	\centering
	\includegraphics[width=0.32\textwidth]{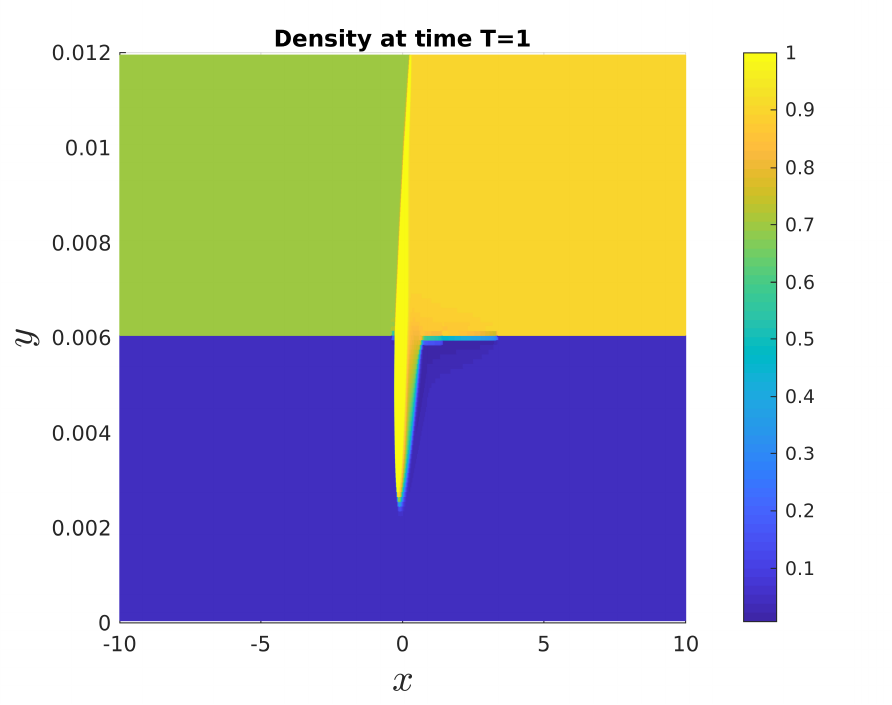}
	\includegraphics[width=0.32\textwidth]{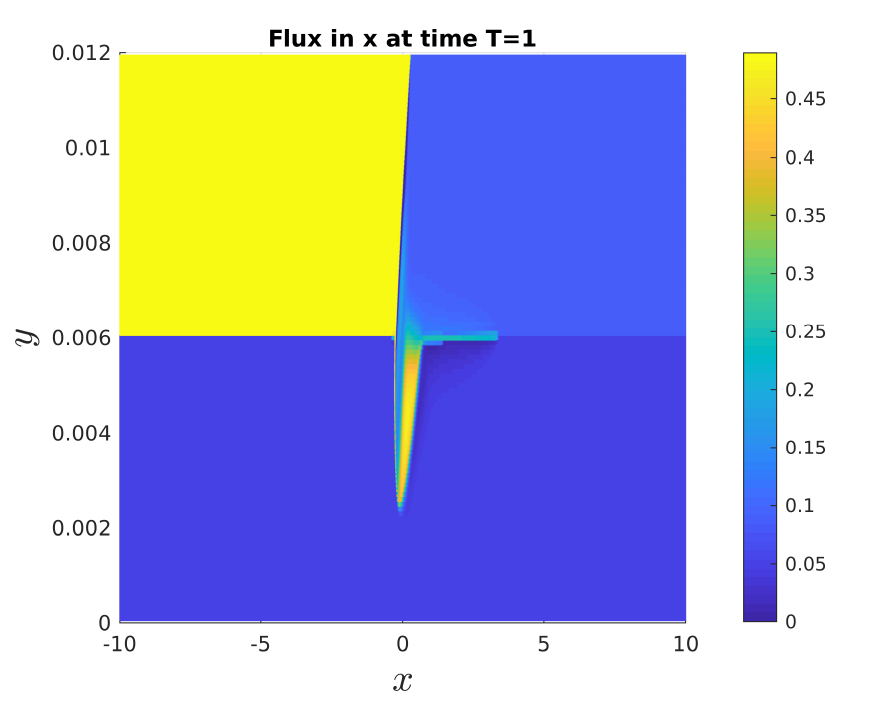}
	\includegraphics[width=0.32\textwidth]{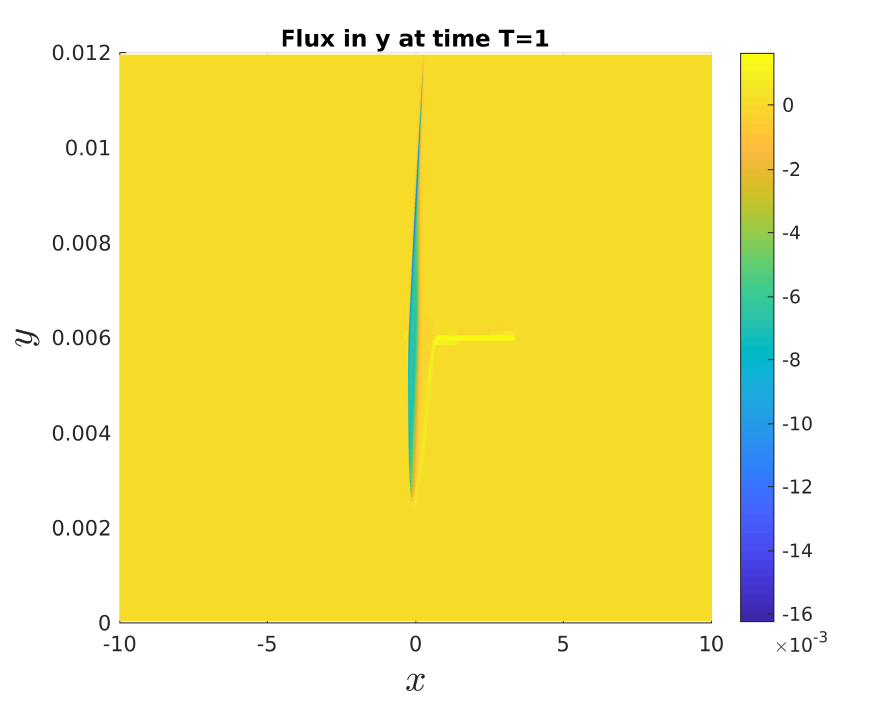}
	\includegraphics[width=0.32\textwidth]{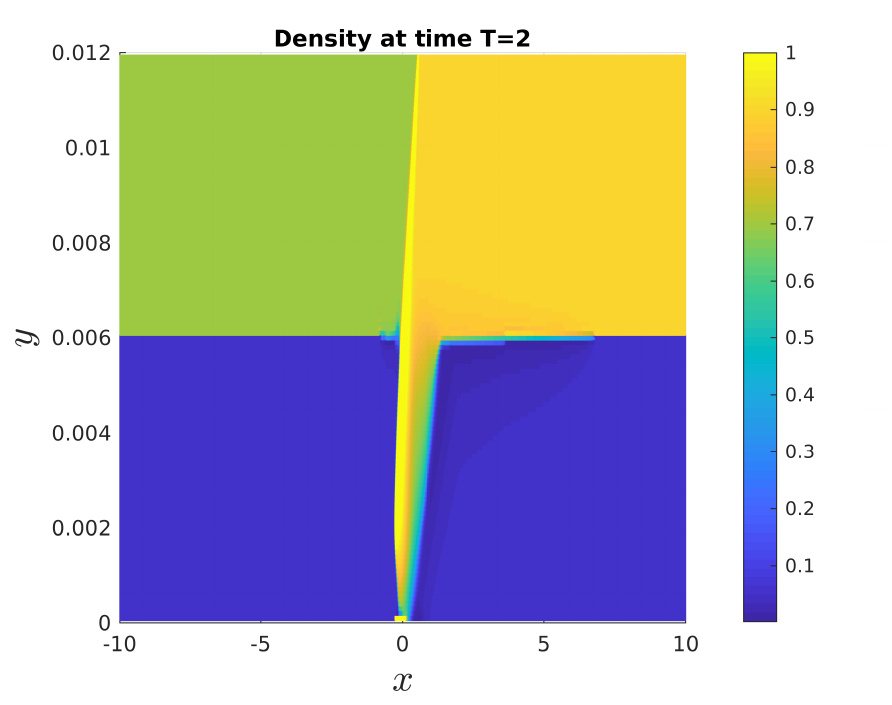}
	\includegraphics[width=0.32\textwidth]{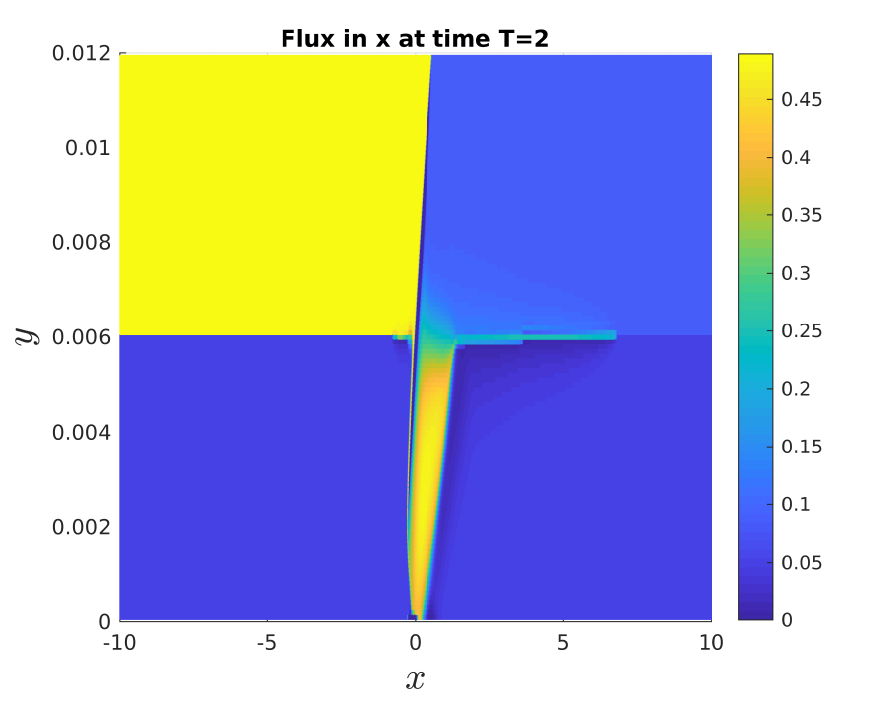}
	\includegraphics[width=0.32\textwidth]{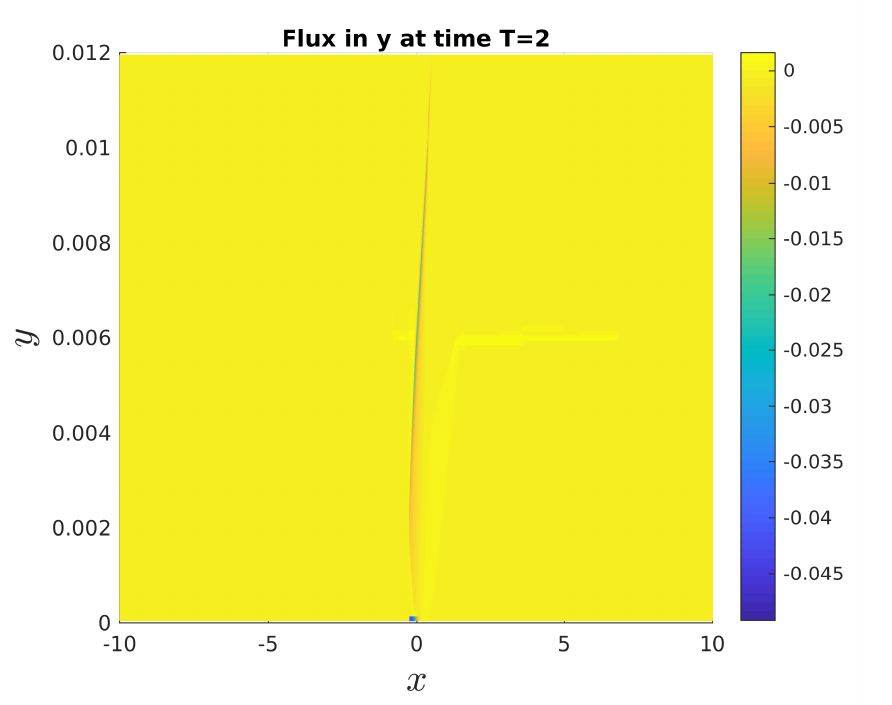}
	\includegraphics[width=0.32\textwidth]{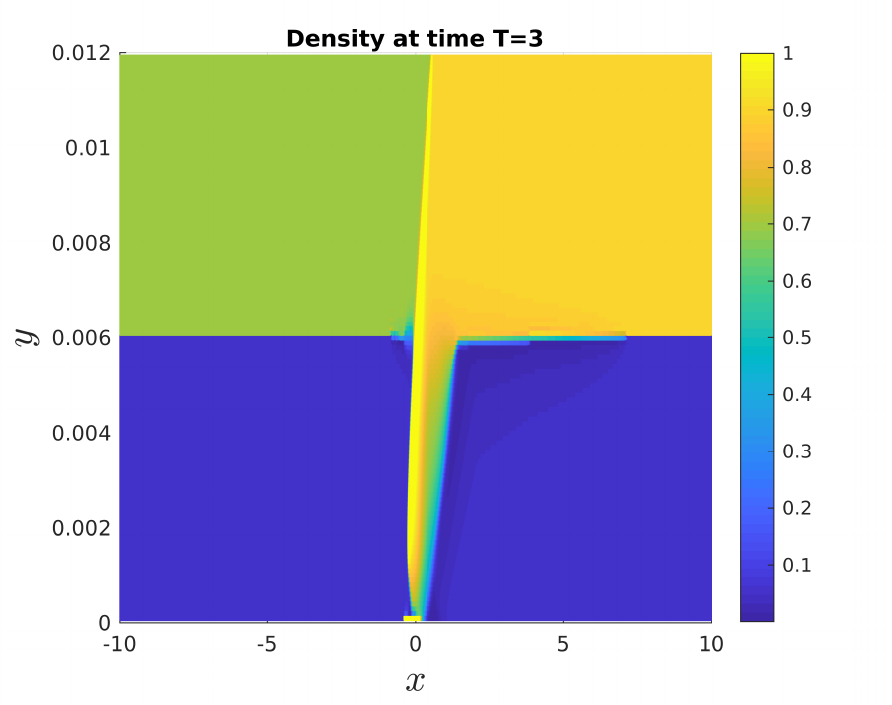}
	\includegraphics[width=0.32\textwidth]{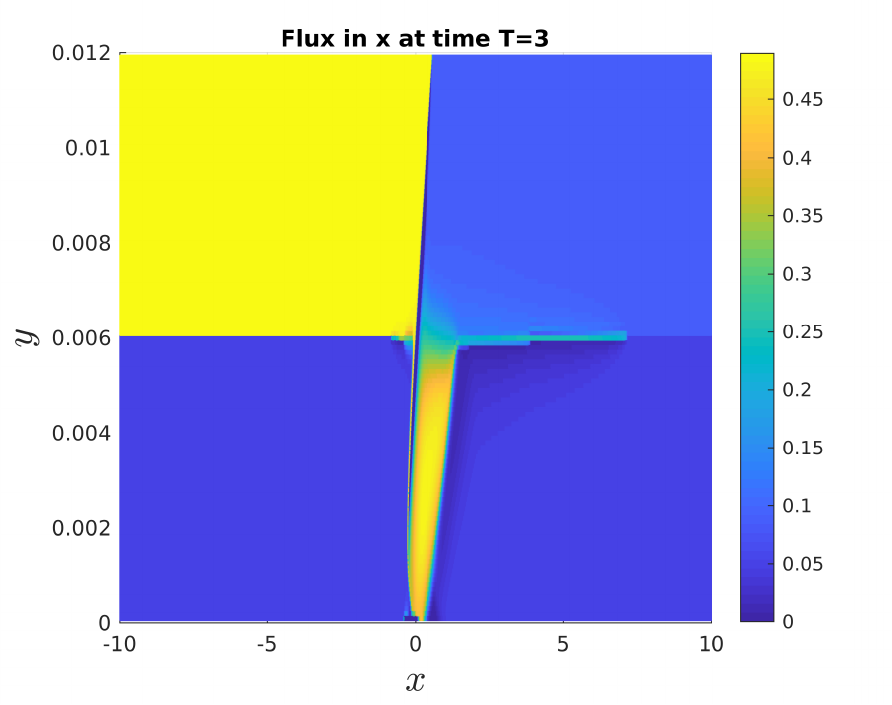}
	\includegraphics[width=0.32\textwidth]{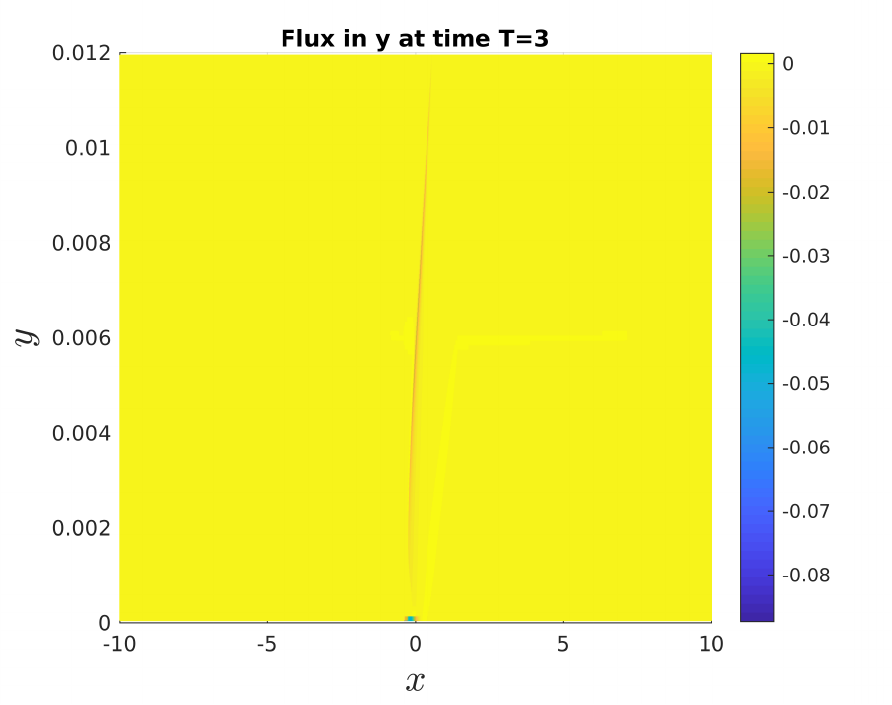}
	\caption{Density $\rho$ (left column), flux $\rho u$ (center column) and flux $\rho v$ (right column) profiles provided by the two-dimensional second-order macroscopic model~\eqref{eq:2DARZEul} at time $T_\text{fin} = 1$ (top row), $T_\text{fin} = 2$ (center row) and $T_\text{fin} = 3$ (bottom row).\label{fig:ToRight}}
\end{figure}

In Figure \ref{fig:ToRight}, we present the density and the flux profiles provided by the 2D ARZ-type model at time $T=1$ (top row), $T=2$ (center row) and $T=3$ (bottom row).

\section{Concluding remarks}\label{sec:conclusion}

This paper introduced a 2D extension of the Aw-Rascle-Zhang~\cite{aw2000SIAP,ZhangMacro} second order model of traffic flow. The proposed model is rather simplistic and it can be viewed as a preliminary step towards multi-lane traffic modeling using 2D second order models. Nonetheless, it enables to capture traffic dynamics caused by lane changing maneuvers. Moreover, it complies with the desirable anisotropic feature of vehicular traffic flow since the wave speeds never exceed the maximum speed of the vehicles. 

Hence, the proposed model opens many perspectives for future research toward several directions. In order to calibrate and test thoroughly the model, enough real data on the traffic macroscopic variables are required. We plan in future work to address this issue again and hope to provide a rigorous validation of the model. Furthermore, the introduction of stochastic features in the lane-changing occurrence, derived from real data, is worth investigating. At the moment, it is clear to us how to relate those quantities to the numerical fluxes in $y-$direction. Finally, a more detailed study on the analytical properties of the model could be provided.

\subsection*{Acknowledgment}
	This work has been supported by HE5386/13-15 and DAAD MIUR project.
    We also thank the ISAC institute at RWTH Aachen, Prof. M. Oeser, MSc. A. Fazekas and MSc. F. Hennecke
    for kindly providing the trajectory data. 

\bibliographystyle{siam}
\bibliography{references}

\begin{thebibliography}{10}

\bibitem{AlNasurThesis}
{\sc S.~Al-nasur}, {\em New Models for Crowd Dynamics and Control}, PhD thesis,
  Faculty of the Virginia Polytechnic Institute and State University, 2006.

\bibitem{AlNasurKachroo}
{\sc S.~Al-nasur and P.~Kachroo}, {\em A microscopic-to-macroscopic crowd
  dynamic model}, in IEEE Conference on Intelligent Transportation Systems,
  2006, pp.~606--611.

\bibitem{aw2002SIAP}
{\sc A.~Aw, A.~Klar, T.~Materne, and M.~Rascle}, {\em Derivation of continuum
  traffic flow models from microscopic follow-the-leader models}, SIAM J. Appl.
  Math., 63 (2002), pp.~259--278.

\bibitem{aw2000SIAP}
{\sc A.~Aw and M.~Rascle}, {\em Resurrection of ``second order'' models of
  traffic flow}, SIAM J. Appl. Math., 60 (2000), pp.~916--938 (electronic).

\bibitem{GoatinChalons}
{\sc C.~Chalons and P.~Goatin}, {\em Transport-equilibrium schemes for
  computing contact discontinuities in traffic flow modeling}, Comm. Math.
  Sci., 5 (2007), pp.~533--551.

\bibitem{ChetverushkinChurbanovaFurmanovTrapeznikova2010}
{\sc B.~N. Chetverushkin, N.~G. Churbanova, I.~R. Furmanov, and M.~A.
  Trapeznikova}, {\em 2{D} micro- and macroscopic models for simulation of
  heterogeneous traffic flows}, in Proceedings of the {ECCOMAS} {CFD} 2010, {V}
  {E}uropean {C}onference on {C}omputational {F}luid {D}ynamics, J.~C.~F.
  Pereira and A.~Sequeira, eds., Lisbon, Portugal, 2010.

\bibitem{ChetverushkinChurbanovaSukhinovaTrapeznikova2008}
{\sc B.~N. Chetverushkin, N.~G. Churbanova, A.~Sukhinova, and M.~Trapeznikova},
  {\em Congested traffic simulation based on a 2{D} hydrodynamical model}, in
  Proceedings of 8th {W}orld {C}ongress on {C}omputational {M}echanics and 5th
  {E}uropean {C}ongress on {C}omputational {M}ethods in {A}pplied {S}cience and
  {E}ngineering, {WCCM8} \& {ECCOMAS} 2008, B.~Schrefler and U.~Perego, eds.,
  Barcelona, Spain, 2008.

\bibitem{Daganzo}
{\sc C.~F. Daganzo}, {\em Requiem for second-order fluid approximation to
  traffic flow}, Transport. Res. B-Meth., 29 (1995), pp.~277--286.

\bibitem{Daganzo2002a}
\leavevmode\vrule height 2pt depth -1.6pt width 23pt, {\em A behavioral theory
  of multi-lane traffic low part i: long homogeneous freeway section},
  Transport. Res. B-Meth., 36 (2002), pp.~131--158.

\bibitem{DiFrancescoFagioliRosini}
{\sc M.~Di~Francesco, S.~Fagioli, and M.~D. Rosini}, {\em Many particle
  approximation of the aw-rascle-zhang second order model for vehicular
  traffic}, Mathematical Biosciences and Engineering, 14 (2017), pp.~127--141.

\bibitem{LawrenceBook}
{\sc L.~C. Evans}, {\em Partial differential equations}, American Mathematical
  Society, 2010.

\bibitem{FanHertySeibold}
{\sc S.~Fan, M.~Herty, and B.~Seibold}, {\em Comparative model accuracy of a
  data-fitted generalized {A}w-{R}ascle-{Z}hang model}, Netw. Heterog. Media.,
  9 (2014), pp.~239--268.

\bibitem{FermoTosin14}
{\sc L.~Fermo and A.~Tosin}, {\em Fundamental diagrams for kinetic equations of
  traffic flow}, Discrete Contin. Dyn. Syst. Ser. S, 7 (2014), pp.~449--462.

\bibitem{Greenbergetal2003}
{\sc J.~M. Greenberg, A.~Klar, and M.~Rascle}, {\em Congestion on multilane
  highways}, SIAM J. Appl. Math., 63 (2003), pp.~813--818.

\bibitem{Helbing}
{\sc D.~Helbing}, {\em Modeling multi-lane traffic flow with queuing effects},
  Physica A: Statistical Mechanics and its Applications, 242 (1997),
  pp.~175--194.

\bibitem{HmVg2017}
{\sc M.~Herty, A.~Fazekas, and G.~Visconti}, {\em A two-dimensional data-driven
  model for trafficflow on highways}, Netw. Heterog. Media, 13 (2018).

\bibitem{hertyillner12}
{\sc M.~Herty and R.~Illner}, {\em Coupling of non-local driving behaviour with
  fundamental diagrams}, Kinet. Relat. Models, 5 (2012).

\bibitem{HmTaVgZm}
{\sc M.~Herty, A.~Tosin, G.~Visconti, and M.~Zanella}, {\em Hybrid stochastic
  kinetic description of two-dimensional traffic dynamics}.
\newblock Preprint: arXiv:1711.02424.

\bibitem{HollandandWoods1997}
{\sc E.~Holland and A.~Woods}, {\em A continuum model for the dispersion of
  traffic on two-lane roads}, Transport. Res. B-Meth., 31 (1997), pp.~473--485.

\bibitem{HoogendoornBovy}
{\sc S.~P. Hoogendoorn and P.~H.~L. Bovy}, {\em Platoon-based multiclass
  modeling of multilane traffic flow}, Networks and Spatial Economics, 1
  (2001), pp.~137--166.

\bibitem{IllnerKlarMaterne}
{\sc R.~Illner, A.~Klar, and T.~Materne}, {\em Vlasov-{F}okker-{P}lanck models
  for multilane traffic flow}, Commun. Math. Sci., 1 (2003), pp.~1--12.

\bibitem{KlarWegener96}
{\sc A.~Klar and R.~Wegener}, {\em A kinetic model for vehicular traffic
  derived from a stochastic microscopic model}, Transport. Theor. Stat., 25
  (1996), pp.~785--798.

\bibitem{klar1999SIAP-1}
\leavevmode\vrule height 2pt depth -1.6pt width 23pt, {\em A hierarchy of
  models for multilane vehicular traffic. {I}. {M}odeling}, SIAM J. Appl.
  Math., 59 (1999), pp.~983--1001 (electronic).

\bibitem{klar1999SIAP-2}
\leavevmode\vrule height 2pt depth -1.6pt width 23pt, {\em A hierarchy of
  models for multilane vehicular traffic. {II}. {N}umerical investigations},
  SIAM J. Appl. Math., 59 (1999), pp.~1002--1011 (electronic).

\bibitem{KlarWegener2000ab}
{\sc A.~Klar and R.~Wegener}, {\em Kinetic derivation of macroscopic
  anticipation models for vehicular traffic}, SIAM J. Appl. Math., 60 (2000),
  pp.~1749--1766 (electronic).

\bibitem{LavalandDaganzo2006}
{\sc J.~Laval and C.~F. Daganzo}, {\em Lane-changing in traffic streams},
  Transport. Res. B-Meth., 40 (2006), pp.~251--264.

\bibitem{LaxBook}
{\sc P.~D. Lax}, {\em Hyperbolic systems of conservation laws and the
  mathematical theory of shock waves}, SIAM, 1973.

\bibitem{Lebacque03}
{\sc J.~P. Lebacque}, {\em Two-phase bounded-acceleration traffic flow model:
  analytical solutions and applications}, Transport. Res. Record, 1852 (2003),
  pp.~220--230.

\bibitem{LebacqueGSOM}
{\sc J.~P. Lebacque and M.~M. Khoshyaran}, {\em A variational formulation for
  higher order macroscopic traffic flow models of the {GSOM} family}, Procedia
  - Social and Behavioral Sciences, 80 (2013), pp.~370--394.

\bibitem{lighthill1955PRSL}
{\sc M.~J. Lighthill and G.~B. Whitham}, {\em On kinematic waves. {II}. {A}
  theory of traffic flow on long crowded roads}, Proc. Roy. Soc. London. Ser.
  A., 229 (1955), pp.~317--345.

\bibitem{Michalopoulosetal1984}
{\sc B.~D. Michalopoulos, P.G. and Y.~Yamauchi}, {\em Multilane traffic
  dynamics: some macroscopic considerations}, Transport. Res. B-Meth., 18
  (1984), pp.~377--395.

\bibitem{PgSmTaVg}
{\sc G.~Puppo, M.~Semplice, A.~Tosin, and G.~Visconti}, {\em Fundamental
  diagrams in traffic flow: the case of heterogeneous kinetic models}, Commun.
  Math. Sci., 14 (2016), pp.~643--669.

\bibitem{PgSmTaVg3}
\leavevmode\vrule height 2pt depth -1.6pt width 23pt, {\em Analysis of a
  multi-population kinetic model for traffic flow}, Commun. Math. Sci., 15
  (2017), pp.~379--412.

\bibitem{PgSmTaVg2}
\leavevmode\vrule height 2pt depth -1.6pt width 23pt, {\em Kinetic models for
  traffic flow resulting in a reduced space of microscopic velocities}, Kinet.
  Relat. Models, 10 (2017), pp.~823--854.

\bibitem{richards1956OR}
{\sc P.~I. Richards}, {\em Shock waves on the highway}, Operations Res., 4
  (1956), pp.~42--51.

\bibitem{seibold2013NHM}
{\sc B.~Seibold, M.~R. Flynn, A.~R. Kasimov, and R.~R. Rosales}, {\em
  Constructing set-valued fundamental diagrams from jamiton solutions in second
  order traffic models}, Netw. Heterog. Media, 8 (2013), pp.~745--772.

\bibitem{SukhinovaTrapeznikovaChetverushkinChurbanova2009}
{\sc A.~Sukhinova, M.~Trapeznikova, B.~N. Chetverushkin, and N.~G. Churbanova},
  {\em Two-{D}imensional {M}acroscopic {M}odel of {T}raffic {F}lows},
  Mathematical Models and Computer Simulations, 1 (2009), pp.~669--676.

\bibitem{VgHmPgTa}
{\sc G.~Visconti, M.~Herty, G.~Puppo, and A.~Tosin}, {\em Multivalued
  fundamental diagrams of traffic flow in the kinetic {F}okker-{P}lanck limit},
  Multiscale Model. Simul., 15 (2017), pp.~1267--1293.

\bibitem{ZhangMacro}
{\sc H.~M. Zhang}, {\em A non-equilibrium traffic model devoid of gas-like
  behavior}, Transport. Res. B-Meth., 36 (2002), pp.~275--290.

\bibitem{Zhuetal}
{\sc H.~B. Zhu, N.~X. Zhang, and W.~J. Wu}, {\em A modified two-lane traffic
  model considering drivers' personality}, Physica A, 428 (2015), pp.~359--367.

\end{thebibliography}

\end{document}